\documentclass{amsart}
\usepackage[all]{xy}
\usepackage{amssymb}
\usepackage{amscd}

\theoremstyle{plain}
\newtheorem{introThm}{Theorem}
\newtheorem{Thm}{Theorem}[section]
\newtheorem{Lem}[Thm]{Lemma}
\newtheorem{Pro}[Thm]{Proposition}
\newtheorem{Cor}[Thm]{Corollary}

\theoremstyle{definition}
\newtheorem{Def}[Thm]{Definition}
\newtheorem{Exa}[Thm]{Example}

\newtheorem{Rem}[Thm]{Remark}
\newtheorem{Ntn}[Thm]{Notation}

\newcommand{\af}{\alpha}

\newcommand{\gm}{\gamma}
\newcommand{\dt}{\delta}
\newcommand{\ep}{\varepsilon}

\newcommand{\io}{\iota}

\newcommand{\ld}{\lambda}

\newcommand{\ph}{\varphi}
\newcommand{\ps}{\psi}

\newcommand{\om}{\omega}

\newcommand{\Gm}{\Gamma}

\newcommand{\Ld}{\Lambda}

\newcommand{\Hnil}{{\mathrm{Hnil}}}
\newcommand{\Dirlim}{\varinjlim}
\newcommand{\Invlim}{\varprojlim}
\newcommand{\ZZ}{{\mathbb{Z}}}
\newcommand{\QQ}{{\mathbb{Q}}}
\newcommand{\CC}{{\mathbb{C}}}

\newcommand{\ts}[1]{{\textstyle{#1}}}
\newcommand{\ds}[1]{{\displaystyle{#1}}}

\newcommand{\GL}{{\mathrm{GL}}}
\newcommand{\TL}{{\mathrm{TL}}}
\newcommand{\Lc}{{\mathrm{Lc}}}
\newcommand{\Lg}{{\mathrm{Lg}}}
\newcommand{\ca}{$C^*$-algebra}

\newcommand{\MM}{{\mathrm{Max}}}
\newcommand{\idc}{_{\circ}}

\begin{document}

\title[Banach Algebras and Rational Homotopy Theory]{Banach
           Algebras and Rational Homotopy Theory}
\author[Lupton]{Gregory Lupton}
\address{Department of Mathematics,
          Cleveland State University,
          Cleveland OH 44115}
\email{G.Lupton@csuohio.edu}
\author[Phillips]{N.~Christopher Phillips}
\address{Department of Mathematics,
         University of Oregon,
         Eugene OR 97403-1222}
\email{ncp@darkwing.uoregon.edu}
\author[Schochet]{Claude L.~Schochet}
\address{Mathematics Department,
         Wayne State University,
         Detroit MI 48202}
\email{claude@math.wayne.edu}
\author[Smith]{Samuel B.~Smith}
\address{Department of Mathematics,
         Saint Joseph's University,
         Philadelphia PA 19131}
\email{smith@sju.edu}
\thanks{Research of the second author
partially supported by NSF grant DMS~0302401.}

\date{16~December 2006}


\begin{abstract}
Let $A$ be a unital commutative Banach
algebra with maximal ideal space $\MM (A).$
We determine the rational
H-type of $\GL_n (A),$ the group of invertible
$n \times n$ matrices with coefficients in $A,$ in terms of the
rational cohomology of $\MM (A).$
We also address an old problem of J.~L.\  Taylor.
Let $\Lc_n (A)$ denote the space of ``last columns'' of $\GL_n (A).$
We construct a natural isomorphism
\[
{\Check{H}}^s (\MM (A); \QQ )
  \cong \pi_{2 n - 1 - s} (\Lc_n (A)) \otimes \QQ
\]
for $n > \frac{1}{2} s + 1$
which shows that the rational cohomology groups of $\MM (A)$
are determined by a topological invariant associated to $A.$
As part of our analysis, we determine
the rational H-type of certain gauge groups
$F (X, G)$ for $G$ a Lie group or,
more generally, a rational H-space.
\end{abstract}

\maketitle

\section{Introduction}

Let $A$ be a unital commutative Banach algebra.
By $\MM (A)$ we denote the set of
maximal ideals of the Banach algebra $A$
topologized with the relative weak$^*$-topology.
Let $\GL_n (A)$ denote
the group of invertible $n \times n$
matrices with coefficients in $A.$
In this paper, we describe the rational homotopy of $\GL_n (A)$
and also that of a second topological space associated to $A.$

The difficulty inherent in describing the
ordinary homotopy theory of the group $\GL_n (A)$
is apparent even in the simplest cases.
For instance, when $A = \CC$ we have
\[
\pi_j (\GL_2 (\CC) ) \cong \pi_j (S^3)
\]
for $j > 1.$
Now $\pi_j (S^3) \neq 0$ for infinitely many
values of $j,$ and these groups are known only for $j < 100$ or so.
Of course this particular problem
becomes tractable after rationalization,
since J.-P.~Serre showed that $\pi_j (S^3)$ is finite for $j > 3.$

In fact,
substantially more can be said about $\GL_n (\CC)$
rationally, using classical results and methods.
For observe that polar decomposition provides a
canonical deformation retraction of
$\GL_n (\CC)$ to the compact Lie group $U_n (\CC).$
By results going back to H.~Hopf, the
rational cohomology algebra of $U_n (\CC)$ is an exterior algebra
generated in
odd degrees $1, 3, \ldots, 2 n - 1.$
{}From these facts,
we deduce a rational homotopy equivalence
\[
\GL_n (\CC) \simeq_{\QQ}
   S^1 \times S^3 \times \cdots \times S^{2 n - 1}.
\]
(See Example~\ref{U_n (C)} below.)
This equivalence actually determines the structure of $\GL_n (\CC)$
as a rational H-space since a homotopy-associative H-space
with only odd rational
homotopy groups is rationally homotopy-abelian.
(See Corollary~\ref{G abelian}.)

Our first main result extends the above analysis
to the general case.
We determine the rational homotopy type of
$\GL_n (A)$ in terms of the rational \v{C}ech cohomology
groups
${\Check{H}}^{*}(\MM (A); \QQ )$
of the space $\MM (A).$
In fact, we determine the full structure of $\GL_n (A)$ as a rational
H-space.
To describe our results, we introduce some notation.

Given a space $X$ with basepoint,
we write $X\idc$ for the path component of the basepoint of $X.$
Recall that if $\pi$ is an abelian group,
then $K (\pi, n)$ is the corresponding Eilenberg-Mac Lane space.
(See Example~\ref{EMspace}.)
If $\pi$ is countable abelian,
by a result of J.~Milnor \cite{Mil2} (also see Milgram~\cite{Milg})
multiplication of loops gives
$K (\pi, n)$ the structure of an abelian
topological group, and this structure is unique up to homotopy.
(Again, see Example~\ref{EMspace}.)
A product $\prod_{j \geq 1} K (\pi_j, j)$ of Eilenberg-Mac Lane spaces
with its standard loop multiplication
consequently also has the structure of an
abelian topological group.
(This multiplication is rarely
unique, even up to homotopy.
See Example~\ref{GEM space}.)

Our first main result is the following.

\begin{introThm}\label{intromain1}
Suppose that $A$ is a unital commutative Banach algebra with
maximal ideal space $\MM (A).$
Let
\[
{\Check{V}}_{h, j} = {\Check{H}}^{2 j - 1 - h} (\MM (A); \QQ ).
\]
Then there is a natural rational H-equivalence
\[
\GL_n (A)\idc \, \simeq_{\QQ} \,
\prod_{j = 1}^{n} \prod_{h = 1}^{2 j - 1} K ({\Check{V}}_{h, j}, h )
\]
where the product of Eilenberg-Mac Lane spaces has the standard
loop multiplication.
Thus $\GL_n (A)\idc$ has the rational H-type of an abelian topological
group.
\end{introThm}

As a consequence of Theorem~\ref{intromain1},
we address an old question of
J.~L.\  Taylor regarding the cohomology of the space $\MM (A).$
In~\cite{Tay}, Taylor raised the
question of characterizing topological invariants
of $\MM (A)$ in terms of invariants of the algebra.
Taylor showed that classical results give a way
to characterize $H^n (\MM (A); \ZZ)$ for $n = 0, 1, 2.$
To be more specific:

\begin{itemize}
\item
$H^0 (\MM (A); \ZZ)$ consists of
formal combinations of idempotents of $A.$
Taylor derives this from
the Shilov Idempotent Theorem.
\item
$H^1 (\MM (A); \ZZ)$ is the
quotient of $\GL (A)$ by the image of the exponential map
$\exp \colon A \to \GL (A).$
Taylor derives this from the Arens-Royden Theorem.
For $A = C (X),$ the algebra of continuous complex valued
functions on a compact Hausdorff space,
and when $X$ is metrizable,
K.~Thomsen observes
(see the introduction to Section~4 of~\cite{T})
that this is a (1934) theorem of Bruschlinsky.
\item
$H^2 (\MM (A); \ZZ)$
is the Picard group ${\mathrm{Pic}} (A)$ of $A$ by a result of Forster.
The Picard group consists of isomorphism
classes of finitely generated projective $A$-modules $P$ which
are invertible in the
sense that there exists some finitely generated
projective module $Q$ such that $P \otimes Q$ is a free module
of rank 1 over $A.$
For a different, and
deeper, version of the Picard group, see~\cite{BGR}.
\end{itemize}

We note that these results are easy in the case $A = C (X)$;
the main point is that they may be
extended to commutative Banach algebras.

In addition, we have the result of J.~Dixmier and A.~Douady~\cite{DD}:
\begin{itemize}
\item
$H^3 (X; \ZZ)$ is isomorphic to
the Brauer group of algebra bundles over $X$ with fibre $\mathcal{K},$
the compact operators on a separable,
infinite-dimensional Hilbert space.
\end{itemize}

As far as we know, there has been no further
progress on describing
$H^n (\MM (A); \ZZ)$ in a concrete fashion.
However, the recent interesting
work of A.~B.\  Thom~\cite{ThA} does shed
additional light on the problem,
albeit in a very abstract context.

As regards $K$-theory, we can add the following results,
which were known to Atiyah, Adams, Swan and others.
\begin{itemize}
\item
$K^0 (\MM (A))$ is isomorphic to
$K_0 (A),$ which is defined (for any ring $A$) as
the quotient of the free abelian group generated by all finitely
generated projective left $A$-modules modulo the isomorphism and
direct sum relations.
\item
$K^1 (\MM (A))$ is isomorphic to
$\GL_{\infty} (A) / \GL_{\infty} (A)\idc,$
the quotient of the infinite general linear
group of $A$ by the connected component of the identity.
\end{itemize}

Now consider the natural action of $\GL_n (A)$ on the right
$A$-module $A^n$ (the direct sum of $n$ copies of $A$).
The space of ``last columns'' $\Lc_n (A)$ is defined
to be the orbit of the last
standard basis vector $e_n \in A^n$ under the action of $\GL_n (A).$
This is a familiar space, studied for example in~\cite{R1}.
With our second main result, we provide an answer to Taylor's question
for the rational (\v{C}ech) cohomology groups of the space $\MM (A).$

\begin{introThm}\label{intromain2}
Let $A$ be a unital commutative Banach algebra.
Then the
rational cohomology groups of $\MM (A)$ are determined
by a topological invariant of $A$ via a natural isomorphism
\[
{\Check{H}}^s (\MM (A); \QQ )
  \cong \pi_{2 n - 1 - s} (\Lc_n (A)) \otimes \QQ
\]
for $n > \frac{1}{2} s + 1.$
\end{introThm}

In an effort to make the paper understandable to both
Banach algebraists and algebraic topologists,
we include more details than would be needed in a paper aimed
at either audience alone.
The organization is as follows.
In Section~\ref{sec:Spaces} we
introduce the topological spaces associated
to a Banach algebra that we will need.
In Section~\ref{sec:Comm Banach} we begin our
examination of the homotopy groups of the
general linear groups of the commutative Banach algebras.
The group $\GL_n (A)$ is known to be homotopy equivalent to the
 space of continuous functions $F (\MM (A), \GL_n (\CC))$
by an important result of Davie~\cite{D}.
We reduce the problem of identifying $\Lc_n (A)$
to one of studying the path component of the constant map in the
function space $F (\MM (A), S^{2 n - 1})$
(Theorem~\ref{commutativeLcn}).
Consequently, we are concerned with studying
function spaces of the form $F (X, G),$ for
$G$ a topological group or, more generally, of the rational homotopy
type of a topological group.
In Section~\ref{sec:rational} we
develop some ideas, mostly from rational homotopy theory,
sufficient to understand completely the rational homotopy type of
a function space of this latter form when $X$ and $G$ are finite
complexes.
The results of this section are developed in greater generality than
strictly necessary, and some are of independent interest for the
homotopy of function spaces (e.g.~Theorem~\ref{Hnil}).
Finally, in Section~\ref{sec:conclusion}, we
prove Theorems~1 and~2 by applying the results developed in
Section~\ref{sec:rational}.

For any topological spaces $X$ and $Y$ we let
$F (X, Y)$ denote the set of continuous functions $X \to Y$ with
the compactly generated topology.
(See Section~IV.3 of~\cite{Wh}.)
This topology agrees with the compact-open topology on
compact sets.
Hence
using this topology rather than the compact-open topology has no
effect on homotopy and singular homology.
If $X$ and $Y$ are based spaces then we let $F_{\bullet} (X, Y)$
denote the subspace of basepoint-preserving maps.
We use $[X, Y]$ to denote the set of based homotopy classes of
based maps from $X$ to $Y.$
When working with a topological group
we take the basepoint to be the identity of the group.
Basepoints are always assumed to be
non-degenerate, that is, the inclusion of
the basepoint into the space is a cofibration,
or equivalently, the pair consisting of the space and its basepoint
is an NDR-pair.
(See page~7 and Section~I.4 of~\cite{Wh}.)
Observe that $f \colon X \to Y$
induces maps $f_* \colon F (W, X) \to F (W, Y)$
and $f^* \colon F (Y, Z) \to F (X, Z)$ given by composition and
pre-composition with $f,$ respectively.
We use this notation
also for the maps induced by $f$ on basepoint-preserving function
spaces and on sets of homotopy classes of maps.

A function space is typically not path connected.
We take the basepoint of $F (X, Y)$ and $F_{\bullet} (X, Y)$
to be the constant map sending $X$ to
the basepoint of to $Y.$
Thus $F (X, Y)\idc$ denotes the space of maps from $X$ to $Y$ which are
freely homotopic to the constant map.
Similarly, $F_{\bullet} (X, Y)\idc$
denotes the space of based, null-homotopic maps.
If $G$ is a topological group then the function space
$F (X, G)$ is one also, with multiplication given by
pointwise multiplication of maps.
In this case, the basepoint is the
constant map carrying $X$ to the identity of $G.$
For any compact space $X$ we let $C (X)$ denote
the \ca\  of continuous complex-valued functions on
$X$ with the operations given pointwise.
We denote by $M_n (A)$ the algebra of $n \times n$
matrices with entries in $A.$
If $A$ is a Banach algebra, then $M_n (A)$ is also
a Banach algebra in an obvious way.

{\textbf{Acknowledgement.}}
We would like to thank Daniel Isaksen and Jim Stasheff for helpful
discussions.

\section{The Spaces}\label{sec:Spaces}

In this section, we let $A$ be any Banach algebra.
While we will
be interested later exclusively in the case $A$ is commutative and
unital, the results of this section can be stated in somewhat greater
generality.
We follow Corach and Larotonda \cite{CL84},~\cite{CL86}
and Rieffel~\cite{R1}
in discussing the spaces associated to $A$ that we will need.
Following Section~3 of~\cite{T},
we handle non-unital algebras by unitizing them
and making slight modifications of the definitions.
In order to preserve functoriality,
we therefore unitize all algebras.
Thus, the definitions of $\GL_n (A)$ and $\Lc_n (A)$ given
below do not agree with those in the introduction.
We will show in Proposition~\ref{P_Unital} below that,
when $A$ is unital, they are nevertheless equivalent.

\begin{Ntn}\label{P_UnitNtn}
If $A$ is a Banach algebra, then we denote by $A^+$ its
unitization $A \oplus \CC,$
with the elements written $a + \ld \cdot 1$ (or just $a + \ld$)
for $a \in A$ and $\ld \in \CC,$
and the multiplication obtained by applying the distributive
law in the obvious way.
(See~Definition~1 in Section~3 of~\cite{BD},
or the end of Section~1.2 of~\cite{Mr}.)
\end{Ntn}

The unitization is again a Banach algebra.
If $A$ is a \ca, then there is a unique choice of norm making
$A^+$ a \ca\  (Theorem~2.1.6 of~\cite{Mr}).

\begin{Rem}\label{P_UnitRmk}
When forming $A^+,$
we add a new identity even if $A$ is already unital.
In this case, $A^+$ is isomorphic to the
Banach algebra direct sum $A \oplus \CC.$
Letting $1_A$ and $1_{A^+}$ denote the identities of $A$ and $A^+,$
the isomorphism is given by
\[
a + \ld \cdot 1_{A^+} \mapsto (a + \ld \cdot 1_A, \ld).
\]
This map is a homeomorphism and if $A$ is a \ca\  it is an isometry.
\end{Rem}

\begin{Def}\label{P_GLDfn}
Let $A$ be a Banach algebra and let $A^+$ denote its unitization
as in Notation~\ref{P_UnitNtn}.
We set
\[
\GL_n (A)
 = \{ a \in M_n (A^+) \colon
      {\mbox{$a$ is invertible and $a - 1 \in M_n (A)$}} \},
\]
with the relative topology from $M_n (A^+).$
\end{Def}

Note that $M_n (A)$ is a subset of $M_n (A^+)$ in an obvious way,
and that, in this definition, $1$~denotes the identity of $M_n (A^+).$

\begin{Rem}\label{P_Action}
Let $A$ be a Banach algebra.
For $a = (a_{j, k}) \in \GL_n (A)$
and $x = (x_1, \ldots, x_n) \in (A^+)^n,$ define $a x$
to be the usual left action of $a \in M_n (A^+),$
the right $A^+$-module endomorphisms of $(A^+)^n$
when elements of $(A^+)^n$ are viewed as column vectors.
This defines a jointly
continuous left action of $\GL_n (A)$ on $(A^+)^n.$
\end{Rem}

\begin{Def}\label{P_SpaceDfn}
Let $A$ be a Banach algebra.
We define
\begin{align*}
\Lg_n (A) & = \{ (x_1, \ldots, x_n) \in (A^+)^n
\colon {\mbox{$x_1, \ldots, x_{n - 1}, x_n - 1 \in A$ and}}   \\
 & \hspace*{4em} {\mbox{there are $y_1, \ldots, y_n \in A^+$ such that
             $\sum_{k = 1}^n y_k x_k = 1$}} \}.
\end{align*}
We define $\Lc_n (A)$ to be the orbit of the last standard
basis vector $e_n = (0, \ldots, 0, 1) \in (A^+)^n$
under the action of $\GL_n (A)$ on $(A^+)^n$ of Remark~\ref{P_Action}.
We let $\gm \colon \GL_n (A) \to \Lc_n (A)$
(or $\gm_A$ if it is necessary to specify the algebra)
be $\gm (a) = a \cdot e_n.$
Finally, we define $\TL_n (A) \subset \GL_n (A)$
to be $\{ a \in GL_n (A) \colon a \cdot e_n = e_n \}.$
All these spaces are given the relative norm topology from
$(A^+)^n$ or $M_n (A^+)$ as appropriate.
\end{Def}

As we will see in Proposition~\ref{P_Unital},
the definitions of $\GL_n (A)$ and $\Lc_n (A)$
simplify to those given in the introduction when $A$ is unital.

\begin{Rem}\label{P_Homog}
One immediately checks that $\TL_n (A)$ is the subgroup of
$\GL_n (A)$ consisting of matrices of the form
\[
\left[
\begin{matrix}
   x  &  0  \\
   c  &  1
\end{matrix}
\right],
\]
where $x \in \GL_{n - 1} (A)$ and $c$ is any row
of elements of $A$ of length $n - 1.$
\end{Rem}

\begin{Lem}\label{P_LcInLg}
Let $A$ be a Banach algebra.
\begin{enumerate}
\item\label{P_L1}
The closed vector subspace
\[
L_n = \{ x = (x_1, \ldots, x_n) \in (A^+)^n \colon x - e_n \in A^n \}
\]
is invariant under the action of $\GL_n (A).$
\item\label{P_L2}
The space $\Lg_n (A)$
is invariant under the action of $\GL_n (A).$
\item\label{P_L3}
We have $\Lc_n (A) \subset \Lg_n (A),$
and $\Lc_n (A)$ is equal to the set of all $n$-tuples of
elements of $A^+$ which occur as the last column of an
element of $\GL_n (A).$
\end{enumerate}
\end{Lem}

\begin{proof}
We prove Part~(\ref{P_L1}).
Let $a \in \GL_n (A)$ and let $x \in L_n.$
Then $a - 1 \in M_n (A),$
from which it easily follows that $(a - 1) x \in A^n.$
Also $x - e_n \in A^n.$
So
\[
a x - e_n = (a - 1) x + (x - e_n) \in A^n.
\]

For Part~(\ref{P_L2}),
let $a \in \GL_n (A)$ and let $x \in \Lg_n (A).$
Then $a \cdot x - e_n \in A^n$ by Part~(\ref{P_L1}).
Choose $y = (y_1, \ldots, y_n) \in (A^+)^n$
such that $y_1 x_1 + \cdots + y_n x_n = 1.$
Regarding $y$ as a row vector,
form the matrix product $z = y \cdot a^{-1} \in (A^+)^n.$
Then one checks that $\sum_{k = 1}^n z_k (a \cdot x)_k = 1.$
So $a \cdot x \in \Lg_n (A).$

We prove Part~(\ref{P_L3}).
To identify $\Lc_n (A)$ as the space of last columns,
simply observe that $a \cdot e_n = x$ if and only if
$(x_1, \ldots, x_n)$ is the last column of $a.$
The first statement
follows from this observation and~(\ref{P_L2}).
\end{proof}

The following theorem is primarily due to Corach and Larotonda
(Theorem~1 of~\cite{CL84}) in the unital case,
with parts due to Thomsen (Section~3 of~\cite{T}) in the
\ca\  case.
We will follow~\cite{T} for most of the proof,
since it uses machinery more familiar in Banach algebras.

\begin{Thm}\label{CR_Serre}
Let $A$ be a Banach algebra.
Then:
\begin{enumerate}
\item\label{P_CR1}  
For each $x \in \Lg_n (A),$
the map $a \mapsto a \cdot x$ is an open mapping
(not surjective) from $\GL_n (A)$ to
\[
L_n = \{ x = (x_1, \ldots, x_n) \in (A^+)^n \colon x - e_n \in A^n \}.
\]
\item\label{P_CR2}  
The map $\gm$ defines a principal locally trivial fibre bundle
with structural group $\TL_n (A).$
\item\label{P_Fibr}  
The sequence
\[
\TL_n (A) \longrightarrow \GL_n (A)
\overset{\gm}{\longrightarrow} \Lc_n (A)
\]
is a Serre fibration,
with base point preserving maps.
\item\label{P_CR4}  
This Serre fibration is natural in $A,$
using the obvious maps (described explicitly in the proof).
\item\label{P_CR5}  
The spaces $\GL_n (A),$ $\TL_n (A),$ and $\Lc_n (A)$
are homeomorphic to open subsets of Banach spaces,
and have the homotopy type of CW~complexes.
\item\label{P_CR6}  
We have $\Lg_n (A)\idc = \Lc_n (A)\idc,$
and the map $\gm$ restricts to an open map of
$\GL_n (A)\idc$ onto $\Lc_n (A)\idc.$
\end{enumerate}
\end{Thm}

\begin{proof}
For \ca s (not necessarily unital),
Part~(\ref{P_CR1}) is a slight strengthening of Lemma~3.3 of~\cite{T}.
We follow the proof of Theorem~8.3 of~\cite{Rf0}.
Let $a \in \GL_n (A)$ and let $x \in \Lg_n (A).$
Set $y = a \cdot x,$ and let $\ep > 0.$
We find $\dt > 0$ such that whenever
$z \in L_n$ and $\| z_k - y_k \| < \dt$ for $1 \leq k \leq n,$
then there is $b \in \GL_n (A)$ such that $\| b - a \| < \ep$
and $b \cdot x = z.$
This will clearly suffice to prove Part~(\ref{P_CR1}).

Now, because $y \in \Lg_n (A)$ (Lemma~\ref{P_LcInLg}(\ref{P_L2})),
there are $r_1, \dots, r_n \in A$ such that
$r_1 y_1 + \dots + r_n y_n = 1.$
Set
\[
M = \max (\| r_1 \|, \ldots, \| r_n \|),
\]
and choose $\dt > 0$ so small that
\[
\dt n^2 M < 1
\hbox{ \ \  and \ \ }
(1 - \dt n^2 M)^{-1} < \frac{\ep}{\| a \|}.
\]
Let $(z_1, \dots, z_n) \in A^n$ satisfy
$\| z_k - y_k \| < \dt$ for all $k.$
Let $c$ be the matrix with $(j, k)$ entry
equal to $(z_j - y_j) r_k.$
Then $c \in M_n (A)$ because $z_j - y_j \in A$ for all $j.$
Also,
\[
c_{j, 1} x_1 + \dots + c_{j, n} x_n
 = (z_j - y_j) r_1 x_1 + \dots + (z_j - y_j) r_n x_n
 = z_j - y_j.
\]
Therefore $(1 + c) y = z,$
whence $b = (1 + c) a$ satisfies $b x = z.$
Finally,
\[
\| c \|
  \leq \sum_{j, k = 1}^n \| (z_j - y_j) r_k \|
  < \dt n^2 M.
\]
Therefore $1 + c$ is invertible in $M_n (A),$
with inverse
\[
(1 + c)^{-1} = 1 - c + c^2 - c^3 + \cdots.
\]
In particular, $(1 + c)^{-1} - 1 \in M_n (A),$ so $c \in \GL_n (A),$
and
\[
\| (1 + c)^{-1} \|
 \leq (1 - \| c \|)^{-1}
 < (1 - \dt n^2 M)^{-1}
 < \frac{\ep}{\| a \|}.
\]
So $b \in \GL_n (A)$ and $\| b - a \| < \ep.$
This proves Part~(\ref{P_CR1}).

For \ca s,
Parts~(\ref{P_CR2}) and~(\ref{P_Fibr}) are in Corollary~3.5 of~\cite{T}
and the discussion after Lemma~3.7 of~\cite{T}
(except for the part about the maps preserving the basepoints,
which is immediate from our conventions on basepoints),
and Part~(\ref{P_CR6}) follows from Lemma~3.7 of~\cite{T}.
For Banach algebras,
the only required change is in the definition of $\Lg_n (A);$
see Definition~\ref{P_SpaceDfn}.

The second statement in
Part~(\ref{P_CR5}) for $\Lc_n (A)$ is Lemma~3.6 of~\cite{T}.
For Banach algebras, we first recall that,
by Corollary~5.5 in Chapter~4 of~\cite{LW},
every open subset of a Banach space has the homotopy type
of a CW~complex.
So it is enough to prove the first statement.
For $\Lc_n (A),$ Part~(\ref{P_CR1}) implies that
$\Lc_n (A)$ is an open subset of the Banach space $L_n$ there.
Also, $\GL_n (A)$ is homeomorphic to the open subset of $A^{n^2}$
consisting of those $(a_{j, k})_{j, k = 1}^n \in A^{n^2}$
such that
\[
\left[
\begin{matrix}
   a_{1, 1} + 1  & a_{1, 2}      & \cdots &
a_{1, n}  \\
   a_{2, 1}      & a_{2, 2} + 1  & \cdots &
a_{2, n}  \\
   \vdots        & \vdots        & \ddots &
\vdots &  \\
   a_{n, 1}      & a_{n, 2}      & \cdots &
a_{n, n} + 1
\end{matrix}
\right]
\]
is invertible in $A^+.$
Similarly, $\TL_n (A)$
is homeomorphic to an open subset of $A^{n (n - 1)}.$

Finally, we consider Part~(\ref{P_CR4}).
Let $A$ and $B$ be Banach algebras, and let $\ph \colon A \to B$
be a continuous homomorphism
(not necessarily unital even if $A$ and $B$ are both unital).
Then there is a continuous unital homomorphism
$\ph^+ \colon A^+ \to B^+$
given by $\ph^+ (a + \ld \cdot 1_A) = \ph (a) + \ld \cdot 1_B.$
One checks immediately that the corresponding map
$\ph^+_n \colon M_n (A^+) \to M_n (B^+)$
sends $\GL_n (A)$ into $\GL_n (B),$
and, using the description given in Remark~\ref{P_Homog},
that it sends $\TL_n (A)$ into $\TL_n (B).$
Similarly, the map
\[
(x_1, \ldots, x_n) \mapsto (\ph^+ (x_1), \ldots, \ph^+ (x_n))
\]
is easily seen to send $\Lg_n (A)$ into $\Lg_n (B)$
and to send $(0, \ldots, 0, 1)$ to $(0, \ldots, 0, 1).$
If $a \in \GL_n (A)$ and $x = (x_1, \ldots, x_n) \in \Lg_n (A),$
then the image of $a \cdot x$ in $\Lg_n (B)$ is equal to
\[
\ph_n^+ (a) \cdot (\ph^+ (x_1), \ldots, \ph^+ (x_n)).
\]
It follows that $\Lc_n (A)$ is sent to $\Lc_n (B),$
and that we have a commutative diagram
\[
\begin{CD}
\TL_n (A) @>>> \GL_n (A) @>{\gm_A}>> \Lc_n (A) \\
@V{\ph_n^+}VV  @VV{\ph_n^+}V  @VVV \\
\TL_n (B) @>>> \GL_n (B) @>>{\gm_B}> \Lc_n (B),
\end{CD}
\]
which is what naturality means.
\end{proof}

The spaces $\GL_n (A)$ and $\Lc_n (A)$ are, in general, not connected.
When we focus on homotopy
theoretic properties of these spaces in subsequent sections,
we will restrict our attention
to the connected components $\GL_n (A)\idc$ and $\Lc_n (A)\idc$
of the basepoints.
The following result justifies this restriction.

\begin{Pro}\label{P_ConnCompSame}
Let $A$ be a Banach algebra.
The connected components of the space
$\GL_n (A)$ are path connected and all homeomorphic.
The same holds for $\Lc_n (A).$
The homeomorphisms can be chosen to be natural
with respect to both a single Banach algebra homomorphism and the
map $\gm$ of the fibration of Theorem~\ref{CR_Serre}(\ref{P_Fibr}).
\end{Pro}

\begin{proof}
The connected components of $\GL_n (A)$ and $\Lc_n (A)$ are open
because these spaces are homeomorphic to open subsets of Banach spaces,
by Theorem~\ref{CR_Serre}(\ref{P_CR5}).
Now note that connected open subsets of Banach spaces
are path connected.

For the second part, let $a \in \GL_n (A).$
Then $b \mapsto a b$ is a homeomorphism from $\GL_n (A)$
to itself which sends $\GL_n (A)\idc$ to a
connected, open, and closed subset of $\GL_n (A),$
which must be the connected component of $\GL_n (A)$ containing $a.$
Similarly, the map $b e_n \mapsto a b e_n,$ for $b \in \GL_n (A),$
is a homeomorphism from $\Lc_n (A)\idc$ to
to the connected component of $\Lc_n (A)$ containing $a e_n.$
By definition, every element of $\Lc_n (A)$ has the form $b e_n$
for some $b,$
so this argument applies to all components.

For naturality, if $\ph \colon A \to B,$ and if $a \in \GL_n (A),$
we take the maps involving $B$ to be multiplication by
$\ph_n^+ (a).$
\end{proof}

As a direct consequence of Theorem~\ref{CR_Serre}, we obtain
a long exact homotopy sequence relating the homotopy
groups of $\GL_{n} (A)$ and $\Lc_{n} (A).$
This sequence
for \ca s is due to Thomsen (Section~3 of~\cite{T}).
Note that, in the statement below,
$\pi_0 (\GL_{n - 1} (A))$ and $\pi_0 (\GL_n (A))$
are possibly nonabelian groups, and that
$\pi_0 (\Lc_n (A))$ is just a set.
That the sequence ends with~$0$ means that
$\pi_0 (\GL_n (A)) \to \pi_0 (\Lc_n (A))$ is surjective.

\begin{Thm}\label{T:les1}
For a Banach algebra $A$ and for each $n > 0$ there
is a long exact sequence
\[
\xymatrix{
\cdots \ar[r] & \pi_k ( \GL_{n - 1} (A) ) \ar[r] & \pi_k ( \GL_{n}
 (A) ) \ar[r] & \pi_k ( \Lc_n (A)) \ar `d[lll] `[llld]_{}[llld] \\
\pi_{k - 1} (\GL_{n - 1} (A)) \ar[r] &
  \ \ \ \   \cdots  \ \ \ \ \ar[r] &
 \pi_1 (\GL_n (A)) \ar[r] & \pi_1 (\Lc_n (A))
\ar `d[lll] `[llld]_{}[llld] \\
\pi_0 (\GL_{n - 1} (A)) \ar[r] & \pi_0 (\GL_n (A)) \ar[r] &
 \pi_0 (\Lc_n (A)) \ar[r] & 0}
\]
which is natural in $A.$
\end{Thm}

\begin{proof}
The long exact sequence in homotopy for the Serre fibration
of Theorem~\ref{CR_Serre}(\ref{P_Fibr}) gives the sequence
of the theorem with $\pi_{*} ( \TL_n (A) )$
in place of $\pi_* ( \GL_{n - 1} (A) ).$
(The map $\pi_0 (\GL_n (A)) \to \pi_0 (\Lc_n (A))$ is
surjective because $\GL_n (A) \to \Lc_n (A)$ is surjective.)
The map $\GL_{n - 1} (A) \to \TL_n (A)$ given by
\[
x \mapsto
\left[
\begin{matrix}
   x  &  0  \\
   0  &  1
\end{matrix}
\right]
\]
is a deformation retraction
(see the description of $\TL_n (A)$ given in Remark~\ref{P_Homog})
and hence induces an isomorphism on homotopy.
This gives the sequence as stated.

Naturality is immediate from naturality of the fibration,
Theorem~\ref{CR_Serre}(\ref{P_CR4}).
\end{proof}

We now turn to the simplifications in the unital case.
Stated informally, if $A$ is unital,
then $\GL_n (A)$ is the invertible group of $M_n (A),$
the space $\Lc_n (A)$ is the set of last columns of
invertible elements of $M_n (A),$ etc.
The proof is immediate, and is omitted.

\begin{Pro}\label{P_Unital}
Let $A$ be a unital Banach algebra,
and identify $A^+$ with the Banach algebra direct sum $A \oplus \CC$
as in Remark~\ref{P_UnitRmk}.
In the obvious way,
further identify $M_n (A^+)$ with $M_n (A) \oplus M_n (\CC)$
and $(A^+)^n$ with $A^n \oplus \CC^n.$
Then:
\begin{enumerate}
\item\label{P_U_GL}
The map $a \mapsto (a, 1)$
defines an isomorphism of topological groups
from the group of invertible elements of $M_n (A)$ to $\GL_n (A).$
\item\label{P_U_Lg}
The map $x \mapsto (x, (0, \ldots, 0, 1))$
defines a homeomorphism from
\[
\left\{ (x_1, \ldots, x_n) \in A^n
\colon {\mbox{There are $y_1, \ldots, y_n \in A$ such that
             $\sum_{k = 1}^n y_k x_k = 1$}} \right\}
\]
to $\Lg_n (A).$
\end{enumerate}
These maps convert the obvious action of the invertible group
of $M_n (A)$ on $\Lg_n (A)$ into the restriction of the action
defined in Remark~\ref{P_Action} to $\GL_n (A).$
In particular, they identify $\Lc_n (A)$ with the orbit of
$(0, \ldots, 0, 1)$ under the action of the invertible group
of $M_n (A),$
and $\TL_n (A)$ with the invertible elements of $M_n (A)$
whose last column is $(0, \ldots, 0, 1).$

These maps are natural for continuous unital homomorphisms
of unital Banach algebras. \qed
\end{Pro}

\section{Commutative Banach algebras}\label{sec:Comm Banach}

Suppose now that $A$ is a commutative unital Banach algebra.
In this case, $\MM (A)$ is a non-empty compact Hausdorff space,
and the Gelfand transform provides an algebra homomorphism
\[
\Gm \colon A \to C (\MM (A)).
\]
See Section~17 of~\cite{BD} or Section~1.3 of~\cite{Mr}.
Recall that if $A = C (X)$ then $\MM (A) = X$
and the Gelfand transform is an isomorphism
(Theorem 2.1.10 of~\cite{Mr}).
In general, $\Gm$ induces a continuous unital homomorphism
\[
\Gm_n \colon A \otimes M_n (\CC)
\longrightarrow M_n (C (\MM (A))) = F (\MM (A), M_n (\CC) )
\]
which restricts to a natural continuous homomorphism
\[
\Gm_n \colon \GL_n (A) \longrightarrow
\GL_n (C (\MM (A))) \, = \, F (\MM (A), \GL_n (\CC) ).
\]
Here we recall that the group operation in the space
$F (\MM (A), \GL_n (\CC))$ is pointwise multiplication of functions.

\begin{Rem}\label{P_GelfIsNat}
For naturality statements, it is important to note that
a unital homomorphism $\ph \colon A \to B$ of commutative unital
Banach algebras induces, in a functorial manner,
a continuous function $\ph^* \colon \MM (B) \to \MM (A).$
If $M \subset B$ is a maximal ideal, then $\ph^* (M) = \ph^{-1} (M).$
If we think of the maximal ideal space as consisting of all
(continuous) unital homomorphisms to $\CC,$
then $\ph^* (\om) = \om \circ \ph.$

It is then immediate that $A \mapsto C ( \MM (A))$ is a
covariant functor.
Thus it makes sense
to say that $\Gm \colon A \to C (\MM (A))$ is natural.
Similarly, $A \mapsto {\Check{H}}^* ( \MM (A); \QQ)$ is a
covariant functor.
\end{Rem}

The following important result of Davie
bridges the Gelfand transform gap.

\begin{Thm}[A.~Davie]\label{Davie1}
Let $A$ be a unital commutative
Banach algebra with maximal ideal space $\MM (A).$
Then the natural homomorphism
\[
\Gm_n \colon \GL_n (A) \longrightarrow F (\MM (A), \GL_n (\CC) ).
\]
induced by the Gelfand transform is a homotopy equivalence.
\end{Thm}

\begin{proof}
This is a special case of Theorem~4.10 of~\cite{D}.
\end{proof}

We note in passing that Davie's work together with very elementary
reasoning will substitute for the Shilov Idempotent Theorem and
for the Arens-Royden Theorem in the
computation of $H^0 (\MM (A); \ZZ)$
and $H^1 (\MM (A); \ZZ)$ respectively.
For our purposes, the significance of this result is that it
identifies $\GL_n (A),$ up to H-equivalence
(as in Definition~\ref{P_Heq}), as a particular gauge group.
We devote the remainder of this section to obtaining a related
identification of $\Lc_n (A)\idc.$

\begin{Thm}\label{Davie2}
Let $A$ be a
unital commutative Banach algebra with maximal ideal space $\MM (A).$
Then the Gelfand transform induces a natural map
\[
L \Gm_n \colon \Lc_n (A)\idc \longrightarrow \Lc_n (C (\MM (A))\idc
\]
which is a homotopy equivalence.
\end{Thm}

\begin{proof}
Set $X = \MM (A).$
Apply Theorem~\ref{T:les1} to $A$ and $C (X),$
and naturality in Theorem~\ref{T:les1} to the Gelfand transform,
obtaining a commutative diagram with long exact rows which ends
with
\[
\begin{CD}
\pi_0 (\GL_{n - 1} (A))
 @>>> \pi_0 (\GL_n (A)) @>>> \pi_0 (\Lc_n (A)) @>>> 0 \\
@VVV  @VVV   @VVV \\
\pi_0 (\GL_{n - 1} (C (X)))
 @>>> \pi_0 (\GL_n (C (X))) @>>> \pi_0 (\Lc_n (C (X))) @>>> 0.
\end{CD}
\]
Theorem~\ref{Davie1} and the Five Lemma imply that
$\pi_k (\Lc_n (A)) \to \pi_k (\Lc_n (C (X)))$
is an isomorphism for all $k \geq 1,$
and a direct argument shows that
$\pi_0 (\Lc_n (A)) \to \pi_0 (\Lc_n (C (X)))$ is a bijection of sets.
Thus, $L \Gm_n$ is a weak homotopy equivalence
(on every connected component, by Proposition~\ref{P_ConnCompSame}).
Since these spaces are homotopy equivalent to CW~complexes
(by Theorem~\ref{CR_Serre}(\ref{P_CR5})),
it follows from Theorem V.3.5 of~\cite{Wh} that
$L \Gm_n$ is a homotopy equivalence.
\end{proof}

We now want to explicitly identify $\Lc_n (C (X)).$

\begin{Lem}\label{ChrisLemma1}
Let $X$ be a compact Hausdorff space.
Identify $\Lg_n (C (X))$ with a subset of $C (X)^n$
as in Proposition~\ref{P_Unital}(\ref{P_U_Lg}),
and further identify $C (X)^n$ with $F (X, \CC^n)$
by sending $(f_1, \ldots, f_n)$ to
$x \mapsto (f_1 (x), \ldots, f_n (x))$
Then the image of $\Lg_n (C (X))$ is exactly $F (X, \CC^n - \{ 0 \}).$
\end{Lem}

\begin{proof}
Let $f \in F (X, \CC^n),$
and write $f (x) = (f_1 (x), \ldots, f_n (x)).$

Suppose $f \in F (X, \CC^n - \{ 0 \}).$
Set
\[
g_k (x) = \frac{\overline{f_k (x)}}{|f_1 (x)|^2 + \cdots + |f_n (x)|^2}.
\]
(The denominator is never zero by hypothesis.)
Then $g_1 f_1 + \cdots + g_n f_n = 1,$ so
$(f_1, \ldots, f_n) \in \Lg_n (C (X)).$

Now suppose that there is $x \in X$ such that $f (x) = 0.$
Then everything in the ideal generated by $f_1, \dots, f_n$
must vanish at $x,$ so $(f_1, \ldots, f_n)$ is not in $\Lg_n (C (X)).$
\end{proof}

\begin{Cor}\label{Chriscor6}
Under the obvious map $F (X, \CC^n) \to C (X)^n,$
the image of the space
\mbox{$F (X, \CC^n - {0})\idc$} is exactly $\Lc_n (C (X))\idc.$
\end{Cor}

\begin{proof}
Combine Lemma~\ref{ChrisLemma1} and Theorem~\ref{CR_Serre}(\ref{P_CR6}).
\end{proof}

The following theorem gives the desired identification
of the space $\Lc_n (A).$

\begin{Thm}\label{commutativeLcn}
Let $A$ be a commutative unital Banach algebra.
Then there is a natural homotopy equivalence
\[
\Lc_n (A)\idc \simeq F (\MM (A), S^{2 n - 1})\idc.
\]
\end{Thm}

\begin{proof}
The Gelfand transform induces a natural homotopy equivalence
\[
\Lc_n (A)\idc \overset{\simeq}{\longrightarrow} \Lc_n (C (\MM (A)))\idc
\]
by Proposition~\ref{Davie2}.
Corollary~\ref{Chriscor6} produces a natural homeomorphism
\[
\Lc_n (C (\MM (A)))\idc \, \approx \, F (\MM (A), \CC^n - \{ 0 \})\idc.
\]
The deformation retraction
\[
\CC^n - \{ 0 \} \longrightarrow S^{2 n - 1}
\]
given by $r (v) = v / \| v \| $ induces a natural homotopy equivalence
\[
F (\MM (A), \CC^n - \{ 0 \})\idc \overset{\simeq}{\longrightarrow}
        F (\MM (A), S^{2 n - 1} )\idc,
\]
and so the composite
\[
(\Lc_n (A))\idc \to \Lc_n (C (\MM (A)))\idc
\to F (\MM (A), \CC ^n - \{ 0 \})\idc \to
    F (\MM (A), S^{2 n - 1})\idc
\]
is the required natural homotopy equivalence.
\end{proof}

\section{H-Spaces, Function Spaces, and Rational
                 Homotopy}\label{sec:rational}

Davie's Theorem (Theorem~\ref{Davie1}) gives a multiplicative
homotopy equivalence between $\GL_n (A)$
and the function space $F (\MM (A), \GL_n (\CC)).$
Further, as mentioned in the introduction, polar decomposition gives a
multiplicative homotopy equivalence
$U_n (\CC) \stackrel{\simeq}{\longrightarrow} \GL_n (\CC).$
In this section, we give a complete description of the rationalization
of the group $F (X, G)$ for $G$ a compact Lie group provided $X$ is a
finite complex.
In fact, as it is necessary
for our analysis, we work in somewhat greater generality, allowing
$G$ to be an H-space with various restrictions.
We also analyze the rationalization of the basepoint-preserving
function spaces $F_{\bullet} (X, G).$
Throughout this section, we thus assume all spaces come
equipped with a fixed, non-degenerate basepoint.

We let $G$ denote a connected H-space (Section~III.4 of~\cite{Wh}).
That is, we assume $G$
comes equipped with a continuous map $\mu \colon G \times G \to G,$
for which there is an element $e \in G$ such that the
restrictions of $\mu$ to $\{ e \} \times G \to G$ and to
$G \times \{ e \} \to G$ are both homotopic to the identity map of $G.$
The map $\mu$ is called the {\emph{multiplication}}, and we
will assume---as for topological groups---that the
identity element $e$ is chosen as the basepoint of $G.$
In what follows, we will write either the pair $(G, \mu)$ or just $G$
depending on whether or not we are concerned with
the particular multiplication for $G$ or just the fact that $G$ is an
H-space.

While we mainly have in mind the case
in which $(G, \mu)$ is a topological group,
our results in this section require only that
$\mu$ be
homotopy-associative and have a homotopy inverse:
a map $\nu \colon G \to G$
such that
$\mu \circ (\nu \times 1) \circ \Delta \colon G \to G$
and $\mu \circ (1 \times \nu) \circ \Delta \colon G \to G$
are both null-homotopic.
That is, we will require
that $(G, \mu)$ be {\emph{group-like}} as defined before
Theorem III.4.14 of~\cite{Wh}.
We begin by recalling some basic facts about group-like spaces
in ordinary homotopy theory.

We will often require that $G$ be a CW~complex.
In this context, recall the following result of I.~James,
Corollaries~1.2 and~1.3 of~\cite{Jam}.

\begin{Pro}[I.~James]\label{James}
Let $(G, \mu)$ be a homotopy-associative H-space.
For any CW~complex $X$ the set $[X, G]$ has the structure of a group.
If $G$ is a CW~complex then $(G, \mu)$ is group-like.
\qed
\end{Pro}

\begin{Def}\label{P_Heq}
If $(G, \mu)$ and $(G', \mu')$ are
H-spaces, then we say a map $F \colon G \to G'$ is an
{\emph{H-map}} if it satisfies
\[
\mu' \circ (F \times F) \sim F \circ \mu \colon G \times G \to G',
\]
where ``$\sim$" denotes homotopy of maps.
An H-map $F \colon G \to G'$ that is a (weak) homotopy
equivalence is called a {\emph{(weak) H-equivalence}}.
We say two H-spaces $G$ and $G'$ are {\emph{H-equivalent}}
if there exists an H-equivalence
$F \colon G \stackrel{\simeq}{\longrightarrow} G'.$
\end{Def}

Next we recall the definition of the homotopy nilpotency of a
group-like space as introduced by Berstein and Ganea~\cite{B-G}.
Suppose $(G, \mu)$ is a group-like space.
Let
\[
\ph_2 \colon G \times G \to G
\]
denote the commutator map
\[
\ph_2 (g_1, g_2) = \mu (\mu (g_1, g_2), \mu (\nu (g_1), \nu (g_2))).
\]
Extend this definition by letting $\ph_1 \colon G \to G$
be the identity map and, for $n > 1,$
defining $\ph_n \colon G^n \to G$ by the rule
\[
\ph_n = \ph_2 \circ (\ph_{n - 1} \times \ph_1).
\]
The following corresponds to Definition~1.7 in~\cite{B-G}.

\begin{Def}\label{P_HNilDfn}
Given a group-like space $(G, \mu),$
the {\emph{homotopical nilpotency}} $\Hnil (G)$ of $(G, \mu)$ is
the least integer $n$ such
that the map $\ph_{n + 1}$ above is null-homotopic.
If $\Hnil (G) = 1$ we say $(G, \mu)$
is {\emph{homotopy-abelian}}
\end{Def}

Finally, recall that the homotopy groups of a homotopy-associative
H-space $(G, \mu)$ come equipped with the Samelson bracket.
(See Section X.5 of~\cite{Wh}.)
This is a bilinear pairing
\[
\langle \, , \, \rangle \colon
 \pi_p (G) \times \pi_q (G) \to \pi_{p + q} (G).
\]

The following result is essentially Theorem~4.6 of~\cite{B-G}.

\begin{Pro}[I.~Berstein and T.~Ganea]\label{nil<=Hnil}
If $(G, \mu)$ is H-equivalent to a loop space, then
$\Hnil (G)$ is an upper bound for the longest nonvanishing
Samelson bracket in $\pi_* (G).$
\end{Pro}

\begin{proof}
By Theorem~4.6 of~\cite{B-G},
$\Hnil (\Omega X)$ is an upper bound for the length of
the longest nonvanishing Whitehead product in $\pi_{*} (X).$
Equivalently, by Theorem X.7.10 of~\cite{Wh},
$\Hnil (\Omega X)$ is an upper bound for the length of
the longest nonvanishing Samelson bracket in $\pi_{*} (\Omega X).$
\end{proof}

\begin{Exa}\label{EMspace}
Recall that an {\emph{Eilenberg-Mac Lane}}
space $K (\pi, n)$ is a CW~complex with only one non-zero homotopy
group, namely $\pi_n (K (\pi, n) ) = \pi.$
(See Section~V.7 of~\cite{Wh}
for the basic theory of Eilenberg-Mac Lane spaces.)
If $\pi$ is abelian and $n \geq 1,$ then $K (\pi, n)$ is homotopy
equivalent to the loop space $\Omega K (\pi, n + 1).$
Thus the Eilenberg-Mac Lane spaces are H-spaces with multiplication
given by the usual product of loops.
By Theorem V.7.13 of~\cite{Wh}
this H-structure is unique
up to homotopy.
If $\pi$ is a countable abelian group then $K (\pi, n)$ is
actually an abelian topological group by
the Corollary to Theorem~3 on page~360 of~\cite{Mil2}
(also see Milgram~\cite{Milg}).
In particular, $\Hnil (K (\pi, n)) = 1.$
\end{Exa}

\begin{Exa} \label{GEM space}
More generally, let $\pi_{1}, \pi_{2}, \ldots$ be any sequence
of countable abelian groups.
Then the product space
$\prod_{j \geq 1} K (\pi_{j}, j)$ also has
the structure of an abelian topological group via the product
multiplication.
We refer to this multiplication as the \emph{standard}
multiplication
for a product of Eilenberg-Mac Lane spaces.
We emphasize that a product of Eilenberg-Mac Lane spaces may admit many
non-H-equivalent H-structures.
(For example, see Proposition~Ia in~\cite{Cur}.)
\end{Exa}

We now turn to the function space $F (X, G).$
If $(G, \mu)$ is an H-space
(respectively, a topological group)
than $F (X, G)$ inherits an H-space structure
(respectively, topological group structure)
from that of $G.$
Specifically, let
\[
\widehat{\mu}
 \colon F (X, G) \times F (X, G) \longrightarrow F (X, G)
\]
be given by
\[
\widehat{\mu} (f_1, f_2) (x) = \mu (f_1 (x), f_2 (x))
\]
for $x \in X$ and $f_1, f_2 \colon X \to G.$
If the maps $f_1$ and $f_2$ preserve basepoints, then so does
$\widehat{\mu} (f_1, f_2).$
Hence $F_{\bullet} (X, G)$ inherits an H-space
structure in the same way.
We recall that, when $X$ is the base space of a principal
$G$-bundle, the group $F (X, G)$ is related to the gauge group
of the bundle and has been studied extensively.
(See, for example,~\cite{AB}.)
We will assume from now on that the
function spaces $F (X, G)$ and $F_{\bullet} (X, G)$
come equipped with the multiplication induced by $G.$
We next collect some basic facts regarding these H-spaces.

\begin{Pro}\label{P_EqOfComps}
Let $X$ be a space and let $G$ be a homotopy-associative H-space.
If either $X$ or $G$ is a CW~complex then
the path components of $F (X, G)$ are naturally homotopy
equivalent to each other and similarly for $F_{\bullet} (X, G).$
\end{Pro}

\begin{proof}
In either case, Proposition~\ref{James} implies that
the set $[X, G]$ of based homotopy classes is a group.
This result is an extension to H-spaces (or to group-like spaces,
as appropriate) of the argument of Proposition~\ref{P_ConnCompSame}.
\end{proof}

\begin{Pro}\label{Milnor}
Let $X$ be a compact space
and let $G$ be a homotopy-associative CW~H-space.
Then $F (X, G)$ and $F_{\bullet} (X, G)$ are H-equivalent
to group-like CW~complexes.
\end{Pro}

\begin{proof}
By Theorem~3 of~\cite{Mil},
the function spaces
have the homotopy type of CW~complexes.
The result follows from Proposition~\ref{James}.
\end{proof}

\begin{Pro}\label{Samequiv5}
Let $X$ be a space and let $G$ be a path connected H-space.
Then there is a natural weak H-equivalence
\[
F (X, G)\idc \simeq_w G \times F_{\bullet} (X, G)\idc,
\]
where the right-hand side has the obvious product H-structure.
If $X$ is a compact metric space and $G$ is a connected CW~complex
this is an H-equivalence that is natural in $X.$
\end{Pro}

\begin{proof}
By assumption the basepoint $x_0 \in X$ is non-degenerate,
so that the inclusion $\{ x_0 \} \to X$ is a cofibration.
This implies that the induced map
\[
F (X, Y) \to F (x_0, Y)
\]
is a fibration for any space $Y.$
(See Theorem I.7.8 of~\cite{Wh}.)
This map is just the evaluation
map $\om \colon F (X, Y) \to Y$ defined by $\om (f) = f (x_0).$
Taking $Y = G$ and restricting, we have a
fibre sequence of path connected spaces
\[
F_{\bullet} (X, G)\idc \stackrel{i}{\longrightarrow} F (X, G)\idc
\stackrel{\om}{\longrightarrow} G.
\]
This fibration has a section $s,$ which
assigns to $g \in G$ the constant function $f (x) = g.$
Therefore the long exact
homotopy sequence breaks into split short exact sequences
\begin{equation}\label{split}
\xymatrix{
   0 \ar[r] & \pi_n (F_{\bullet} (X, G)\idc)
   \ar[r]^{i_{*}} &
   \pi_n (F (X, G)\idc) \ar[r]^-{\om_{*}} &
   \pi_n (G) \ar@/_1.5pc/[l]_-{s_{*}} \ar[r] &
   0
}
\end{equation}
for each $n.$
Define an H-map
$\ph \colon G \times F_{\bullet} (X, G)\idc \to F (X, G)\idc$
by $\ph = \widehat{\mu} \circ (s \times i),$ where
$\widehat{\mu}$ is the multiplication on $F (X, G)\idc$
described before Proposition~\ref{P_EqOfComps}.
It is straightforward to check that $\ph$ induces an
isomorphism on homotopy groups.
Finally, if $X$ and $G$ satisfy the further hypotheses above then,
by Proposition~\ref{Milnor},
the function spaces in question
have the homotopy type of CW~complexes.
So a weak homotopy equivalence is a homotopy equivalence,
by Theorem V.3.5 of~\cite{Wh}.
\end{proof}

Regarding the homotopical nilpotency of $F (X, G)\idc,$
we have the following result.

\begin{Thm}\label{Hnil}
Let $X$ be a space and $G$ a connected topological group, or,
alternately
let $X$ be a compact space and $G$ a connected
homotopy-associative CW~H-space.
Then
\[
\Hnil (F (X, G)\idc) = \Hnil (G).
\]
\end{Thm}

\begin{proof}
With either hypothesis,
$F (X, G)\idc$ is group-like and so its homotopical
nilpotency is well-defined.
(In the second case, see Proposition~\ref{Milnor}.)
As in the proof of Proposition~\ref{Samequiv5}, we have a section
$s \colon G \to F (X, G)\idc$ of
the evaluation fibration $\om \colon F (X, G)\idc \to G.$
It is easy to see that both $\om$ and $s$ are
H-maps---recall that
the multiplication on $F (X, G)\idc$ is
induced from that of $G.$
Let
\[
\widehat{\ph}_{n} \colon F (X, G)\idc^{n} \to F (X, G)\idc
\]
denote the commutator map for $F (X, G)\idc.$
Then we have
\[
\widehat{\ph}_{n} \circ s^{n} = s \circ \ph_n
  \colon G^{n} \to F (X, G)\idc.
\]
The inequality $\Hnil (F (X, G)\idc) \geq \Hnil (G)$ follows.
On the other hand, suppose $\Hnil (G) = n$ so
that $\ph_{n + 1} \colon G^{n + 1} \to G$ is null-homotopic.
We can identify
\[
F (X, G)\idc^{n + 1} \cong F (X, G^{n + 1})\idc.
\]
The map $\widehat{\ph}_{n + 1}$ is then adjoint,
via the exponential law, to the composition
\[
\xymatrix{X \times F (X, G^{n + 1})\idc
\ar[r]^-{\ep}
& G^{n + 1} \ar[r]^-{\ph_{n + 1}} & G }
\]
where $\ep (x, f) = f (x)$ is the
generalized evaluation map.
Since this composition is null-homotopic,
$\widehat{\ph}_{n + 1}$ is null-homotopic.
\end{proof}

We now consider the H-space $F (X, G)\idc$ in the context of
rational homotopy theory.
To do so, however,
we will need to make some further restrictions on $X$ and $G.$
In particular, for the remainder of this section, all spaces
(except function spaces)
will be taken to be CW~complexes.

We recall the basic facts about the rationalization of groups and
spaces as developed in~\cite{HMR}.
A connected CW~complex $X$ is
{\emph{nilpotent}} (Definition II.2.1 of~\cite{HMR})
if $\pi_1 (X)$ is a nilpotent group
and $\pi_1 (X)$ acts nilpotently (Section I.4 of~\cite{HMR})
on each $\pi_j (X)$ for $j \geq 2.$
A nilpotent space $X$ is a {\emph{rational space}} if its homotopy
groups $\pi_j (X)$ are $\QQ$-vector spaces for each $j \geq 1.$
(This is the case $P = \varnothing$ of Definition II.3.1 of~\cite{HMR}.)
If $\pi$ is a nilpotent group,
then we denote by $\pi_{\QQ}$ its rationalization
in the sense of Definition I.1.1 of~\cite{HMR}
(that is, the case $P = \varnothing$ there).
If $\pi$ is abelian,
such as for $\pi = \pi_j (X)$ with $j \geq 2$
and also for $\pi_1$ of an H-space,
then $\pi_{\QQ}$ is just $\pi \otimes \QQ.$
If $X$ is a nilpotent CW~complex, then
a function $f \colon X \to Y$ is a {\emph{rationalization}} if $Y$ is
rational and if $f$ induces an isomorphism
\[
f_{\sharp} \colon \pi_j (X)_{\QQ} \stackrel{\cong}{\longrightarrow}
\pi_j (Y)
\]
for each $j \geq 1.$
This is not the same statement as Definition II.3.1 of~\cite{HMR},
but is equivalent by Theorem~3B in Section~II.3 of~\cite{HMR}.

With this terminology, we have the following theorem.

\begin{Thm}[Hilton, Mislin, and Roitberg]\label{rationalizenilpotent}
Every nilpotent CW~complex $X$ has a
rationalization $e \colon X \to X_{\QQ}$
where $X_{\QQ}$ is a CW~complex.
The space $X_{\QQ}$ is unique up to homotopy equivalence.
\end{Thm}

\begin{proof}
Existence is Theorem~3B in Section~II.3 of~\cite{HMR}.
For uniqueness,
suppose $f \colon X \to Y$ is some other rationalization.
Following Definition II.3.1 of~\cite{HMR},
$f$ induces a bijection
\[
f^* \colon [Y, X_{\QQ}] \to [X, X_{\QQ}]
\]
Thus there exists $h \colon Y \to X_{\QQ}$
such that $h \circ f \simeq e.$
Using Theorem~3B in Section~II.3 of~\cite{HMR},
we see that $h$ induces an isomorphism of rational homotopy groups
since $f$ and $e$ do.
Since $Y$ and $X_{\QQ}$ are rational spaces, $h$ is a weak
equivalence.
Finally, since $Y$ and $X_{\QQ}$ are CW~complexes,
Theorem V.3.5 of~\cite{Wh} implies that $h$ is a homotopy equivalence.
\end{proof}

We use the notation $X \simeq_{\QQ} Y$ to mean that $X$ and $Y$ are
{\emph{rationally equivalent}} spaces, that
is, that $X_{\QQ}$ and $Y_{\QQ}$ are homotopy equivalent.
We give some examples of rationalization maps which will figure
prominently below.

\begin{Exa}\label{odd-sphere}
Consider an odd-dimensional sphere $S^{2 n - 1}.$
The groups $\pi_j (S^{2 n - 1})$
are finite except in the single degree $j = 2 n - 1.$
(See Theorem 9.7.7 of~\cite{Sp}.)
So the rationalization of $S^{2 n - 1}$ is a map
\[
e \colon S^{2 n - 1} \to K (\QQ, 2 n - 1)
\]
corresponding to a nontrivial class in $H^{2 n - 1} (S^{2 n - 1}; \QQ).$
Thus, we have
\[
S^{2 n - 1} \simeq_{\QQ} K (\ZZ, 2 n - 1)
   \simeq_{\QQ} K (\QQ, 2 n - 1),
\]
and so, by Example~\ref{EMspace},
an odd-dimensional sphere has the rational homotopy
type of an abelian topological group.
We note in passing that even-dimensional
spheres are slightly more complicated,
having two nontrivial rational homotopy groups.
\end{Exa}

There are many more spaces, like $S^{2 n - 1},$
whose rationalizations are Eilenberg-Mac Lane spaces
or products of them.

\begin{Def}\label{P_RatHD}
Let $X$ be a connected nilpotent CW~complex.
We say $X$ is a {\emph{rational H-space}} if $X_{\QQ}$ is an H-space.
\end{Def}

The following lemma and its corollary are well-known.
(See, for instance,~\cite{Sch}.)
We include the statements and proofs for completeness.

\begin{Lem}\label{G decomposed}
Let $G$ be a connected CW~H-space.
Then $G$ is nilpotent,
and the rationalization of $G$ may be written in the form
\begin{equation}\label{G equiv}
e \colon G \to \prod_{j \geq 1} K (\pi_j (G) \otimes \QQ, j).
\end{equation}
In particular, the rational homotopy groups
of $G$ correspond to a space of algebra generators of $H^{*} (G; \QQ).$
This rationalization is natural in the sense that a
map of H-spaces $f \colon G_{1} \to G_{2}$ gives rise to a
homotopy-commutative diagram
in which the map on Eilenberg-Mac Lane spaces
is induced by
$f_{\sharp} \otimes 1 \colon
\pi_{*} (G_{1}) \otimes \QQ \to \pi_{*} (G_{2}) \otimes \QQ.$
\end{Lem}

\begin{proof}
The fundamental group of $G$ is abelian and acts trivially,
that is, $G$ is a {\emph{simple}} space.
(See Theorem 7.3.9, and the preceding discussion, in~\cite{Sp}.)
So $G$ is nilpotent.
By a classical theorem of H.~Hopf
(Corollary III.8.12 of~\cite{Wh}),
the rational cohomology algebra $H^* (G; \QQ)$
is the tensor product of an exterior algebra on odd generators
with a polynomial algebra on even generators,
that is, $H^{*} (G; \QQ)$ is a free commutative graded algebra.
(See Example~6 in Section~3 of~\cite{FHT}.)
Write $n_{1}, n_{2}, \ldots$ for the degrees of these generators.
We then obtain a map $e \colon G \to \prod_{j} K (\QQ, n_{j})$
such that, on rational cohomology,
the image of the fundamental class of
$K (\QQ, n_{j})$ is the corresponding generator of $H^{*} (G; \QQ).$
Since $H^{*} (G; \QQ)$ is free over $\QQ,$
the map $e$ induces an isomorphism on rational cohomology.
By Theorem~3B in Section~II.3 of~\cite{HMR} again,
$e$ is a rational equivalence.
\end{proof}

Thus we see that it is sufficient to determine the rational
homotopy groups of an H-space, or more generally of a rational
H-space, in order to identify its rational homotopy type.

\begin{Cor}\label{G_Dec_Cor}
Let $X$ be a rational H-space (Definition~\ref{P_RatHD}).
Then there exists a rational equivalence of the form
\[
X \simeq_{\QQ} \prod_{j\geq 1} K (\pi_{j} (X) \otimes \QQ, j),
\]
which is natural in the same sense as in Lemma~\ref{G decomposed}.
\end{Cor}

\begin{proof}
Apply Lemma~\ref{G decomposed} to $X_{\QQ},$
which is a CW~complex by Theorem~\ref{rationalizenilpotent},
and is still connected,
using the fact that the rationalization of a rational space is a
homotopy equivalence.
\end{proof}

\begin{Exa} \label{U_n (C)}
We specialize Lemma~\ref{G decomposed} to the Lie group $U_{n} (\CC).$
The cohomology $H^* (U_n (\CC); \QQ)$ is the exterior
algebra on one generator in each odd degree from~$1$ through $2 n - 1.$
(See page~412 of~\cite{Br}.)
Thus, using Example~\ref{odd-sphere},
\[
U_n (\CC) \simeq_{\QQ} S^1 \times S^3 \times \cdots
 \times S^{2 n - 1} \simeq_{\QQ} \prod_{j = 1}^{n} K (\QQ, 2j-1),
\]
as mentioned in the introduction.
Alternately,
we may compute the rational homotopy of $U_n (\CC)$
from the long exact sequences of the fibrations
$U_{k - 1} (\CC) \longrightarrow U_k (\CC) \longrightarrow S^{2 k - 1}$
for $2 \leq k \leq n$
(used in the analysis in Section~IV.10 of~\cite{Wh}).
In particular, by IV.10.17 of~\cite{Wh},
the map $U_{n - 1} (\CC) \to U_n (\CC)$
is an isomorphism on $\pi_k$ for $k \leq 2 n - 3.$
\end{Exa}

\begin{Rem}\label{remark:H}
The preceding examples display various possible
situations which can occur in the context of H-spaces and
rationalization.
As we show below (Corollaries~\ref{P_Loop} and~\ref{G abelian}),
the rational equivalence
$U_n (\CC) \simeq_{\QQ} \prod_{j = 1}^{n} K (\QQ, 2j-1)$
can be taken to be multiplicative when the
product of Eilenberg-Mac Lane spaces has the standard multiplication.
On the other hand, for a general H-space $G,$ the map
$G \to \prod_{j \geq 1} K (\pi_j (G) \otimes \QQ, j)$
of Lemma~\ref{G decomposed} will rarely be  multiplicative or
even an $H$-map,
even up to rational equivalence, if the product of Eilenberg-Mac Lane
spaces has the standard multiplication.
(See Example~\ref{notratH}.)
Finally, the space $S^{2 n - 1}$ is not
an H-space for $n \neq 1, 2, 4$ (by J.~F.\  Adams~\cite{Ad})
and so the question
of the multiplicativity of $S^{2 n - 1} \to K (\QQ, 2 n - 1)$
in these cases is moot.
\end{Rem}

Turning to the rationalization of function spaces, the
results we need depend on three basic developments.
First, as has already been used, Theorem~3 of~\cite{Mil}
shows that the path components of
$F (X, Y)$ have the homotopy type of CW~complexes provided $X$ is
a compact metric space.
(Recall we assume $Y$ has the homotopy
type of a CW complex.)
The second ingredient is Theorem II.3.11 of~\cite{HMR}
and the discussion that follows it.
Nilpotence of the spaces is Corollary II.2.6 of~\cite{HMR}.

\begin{Thm}[Hilton, Mislin, and Roitberg]\label{HMRlocal}
Suppose that $X$ is a finite CW~complex
and $Y$ is a nilpotent space with
rationalization $e \colon Y \to Y_{\QQ}.$
Then the components of $F(X, Y)$ and $F_{\bullet}(X, Y)$ are all
nilpotent spaces and the induced maps
\[
e_* \colon F_{\bullet} (X, Y)\idc \longrightarrow
F_{\bullet} (X, Y_{\QQ})\idc
\hbox{\, \and \, }
e_* \colon F (X, Y)\idc \longrightarrow F (X, Y_{\QQ})\idc
\]
given by composition with $e$ are rationalizations.
\qed
\end{Thm}

The third ingredient is an early result of R.~Thom
(Theorem~2 of~\cite{Th}).

\begin{Thm}\label{oldThom}
(R.~Thom)
Let $X$ be a Hausdorff topological space.
Then there is a natural isomorphism
\[
\pi_j (F (X, K (\pi, n))\idc) \cong H^{n - j} (X; \pi)
\]
for $j \geq 1.$
\qed
\end{Thm}

If $X$ is not a finite CW~complex then one
must specify which cohomology theory is being used:
Thom is using singular cohomology.
This is not the best choice for compact
spaces; \v{C}ech theory is better, and in general the natural
map from \v{C}ech theory to singular theory
is neither injective nor surjective.
We use Thom's result only for $X$ a finite complex.

Combining the facts above, we obtain the following result.

\begin{Thm}\label{null-maps}
Let $X$ be a finite CW~complex and let $G$
be a connected CW~H-space or, more generally, a rational H-space.
Let $V_j$ and $\widetilde{V}_j$ be the rational vector spaces
\[
V_j
 \, = \, \bigoplus_{l \geq j}
      H^{l - j} (X; \pi_{l} (G) \otimes \QQ)
{\hbox{ \ \ and \ \ }}
\widetilde{V}_j \, = \,
   \bigoplus_{l \geq j}
      \widetilde{H}^{l - j} (X; \pi_{l} (G) \otimes \QQ).
\]
Then there are natural rational equivalences
\[
F (X, G)\idc \, \simeq_{\QQ} \, \prod_{j \geq 1} K (V_j, j)
\hbox{\ \ and \ \ }
F_{\bullet} (X, G)\idc \, \simeq_{\QQ} \,
   \prod_{j \geq 1} K (\widetilde{V}_j, j).
\]
\end{Thm}

\begin{proof}
Since $G$ is a rational H-space, we apply Corollary~\ref{G_Dec_Cor}
to write its
rationalization as a product of Eilenberg-Mac Lane spaces
\[
G_{\QQ} \, \simeq \, \prod_{l \geq 1} K (\pi_l (G) \otimes {\QQ}, l).
\]
By Theorem~\ref{HMRlocal}, we have rational equivalences
\[
F (X, G)\idc \simeq_{\QQ} F (X, G_{\QQ})\idc
{\hbox{ \, and \, }}
F_{\bullet} (X, G)\idc \simeq_{\QQ} F_{\bullet} (X, G_{\QQ})\idc.
\]
Now $F (X, G_{\QQ})\idc$ and
$F_{\bullet} (X, G_{\QQ})\idc$ are H-spaces.
Thus, by Lemma~\ref{G decomposed}, to prove the theorem
it suffices to compute the (rational) homotopy groups of these spaces.

Using the
standard homeomorphism of function spaces
$F (X, \prod_l Y_l) \approx \prod_l F (X, Y_l),$
we obtain
\[
F (X, G_{\QQ})\idc
 \simeq F \left( X, \ts{\ds{\prod}_{l \geq 1}}
           K (\pi_l (G) \otimes {\QQ}, l) \right)\idc
 \approx \prod_{l \geq 1}
           F \left( X, K (\pi_l (G) \otimes {\QQ}, l) \right)\idc.
\]
Thom's result~\ref{oldThom} now gives
\[
\pi_{j}\left( F \big( X, K (\pi_l (G) \otimes {\QQ}, l) \big)\idc
\right)
\cong
H^{l - j} (X; \pi_l (G) \otimes {\QQ}).
\]
Thus we see $\pi_{j} (F (X, G_{\QQ})\idc) \cong V_{j},$ as needed.
Next, using this isomorphism
and the split short exact sequence~(\ref{split})
we compute
$\pi_{j} (F_{\bullet} (X, G_{\QQ})\idc) \cong \widetilde{V}_{j},$
which completes the proof.
\end{proof}

The ingredients of Theorem~\ref{null-maps}
have been combined before in the same way,
for specific $X$ and $G.$
(See, for instance, the computations in Section~2 of~\cite{AB}.)
While Theorem~\ref{null-maps}
determines the full rational homotopy type of the function spaces of
interest to us,
these function spaces are H-spaces and so it is natural to ask
whether we can determine their structure as such,
at least up to rational
H-equivalence.
We take up this question for the remainder of the section.

Following the discussion before Theorem II.1.8 of~\cite{HMR},
if $(G, \mu)$ is an H-space then the rationalization $G_{\QQ}$
admits a multiplication $\mu_{\QQ}$ making the rationalization map
$e \colon G \to G_{\QQ}$ an H-map.
Moreover, this H-space structure
on $G_{\QQ}$ is uniquely determined by that of $G$ in the sense that if
$(\widetilde{G}_{\QQ}, \widetilde{\mu}_{\QQ})$
represents a second rationalization of $(G, \mu)$ then
$(G_{\QQ}, \mu_{\QQ})$ and
$(\widetilde{G}_{\QQ}, \widetilde{\mu}_{\QQ})$
are H-equivalent.
It is easy to see that if $(G, \mu)$
is group-like, then so is $(G_{\QQ}, \mu_{\QQ}).$

\begin{Def} \label{RatHequiv}
We say two CW~H-spaces $(G, \mu)$ and $(G', \mu')$
are {\emph{rationally H-equivalent}} if
their rationalizations $(G_{\QQ}, \mu_{\QQ})$ and
$(G'_{\QQ}, \mu'_{\QQ})$ are H-equivalent.
\end{Def}

\begin{Exa} \label{notratH}
The rationalization $K (\ZZ, n) \to K (\QQ, n)$ is clearly a rational
H-equivalence.
On the other hand, the rationalization
$\Omega S^{2} \to K (\QQ, 1) \times K (\QQ, 2)$
is not a rational H-equivalence if both sides
have the standard loop multiplications.
To see that these spaces are not H-equivalent, we argue as follows.
Let $\io \in \pi_2 (S^2)$ be the homotopy class of the identity map.
As is well-known,
$[\iota, \iota] \in \pi_3(S^2)$ is of infinite order
(see Theorem XI.2.5 of~\cite{Wh}), yielding a
non-zero Whitehead bracket in $\pi_* (S^2_{\QQ}).$
Under the identification of
Whitehead brackets in $\pi_* (S^2_{\QQ})$ with Samelson brackets in
$\pi_* (\Omega S^2_{\QQ})$ (Theorem X.7.10 of \cite{Wh}), there is a
corresponding non-zero Samelson bracket in
$\pi_* (\Omega S^2_{\QQ}).$
By Proposition~\ref{nil<=Hnil}
we have $\Hnil(\Omega S^2_{\QQ}) \geq 2.$
On the other hand, $K (\QQ,1) \times K (\QQ,2)$ is homotopy-abelian.
\end{Exa}

\begin{Def}\label{ratHnil}
Given a connected CW~homotopy-associative H-space $(G, \mu),$
the {\emph{rational homotopical nilpotency}} $\Hnil_{\QQ} (G)$
is the homotopical nilpotency of the
rationalization $(G_{\QQ}, \mu_{\QQ})$ of $(G, \mu),$ that is,
$\Hnil_{\QQ} (G) = \Hnil (G_{\QQ}).$
We say $(G, \mu )$ is
{\emph{rationally homotopy-abelian}} if $\Hnil_{\QQ} (G) = 1.$
\end{Def}

The Samelson bracket in $\pi_* (G)$ induces a bracket in the
rationalization giving $\pi_* (G) \otimes \QQ$ the structure of
a {\emph{graded Lie algebra}} as
in the following definition.
(See the beginning of Section~21 of~\cite{FHT}.)

\begin{Def}\label{gla}
A positively graded vector space $L$ over $\QQ$ is a
{\emph{graded Lie algebra}} if
$L$ comes equipped with a bilinear, degree zero pairing
$\langle \, , \, \rangle$
satisfying
\begin{enumerate}
\item[(1)] Anti-symmetry:
$\langle \alpha, \beta \rangle =
 - (-1)^{\deg (\alpha) \deg(\beta)} \langle \beta, \alpha \rangle$
\item[(2)] Jacobi identity:
$\langle \alpha, \langle \beta, \gamma \rangle \rangle
 = \langle \langle \alpha, \beta \rangle, \gamma\rangle
 + (-1)^{\deg (\alpha) \deg(\beta)}
  \langle \beta, \langle \alpha, \gamma \rangle \rangle.$
\end{enumerate}
\end{Def}

By Theorems X.5.1 and X.5.4 of~\cite{Wh},
after tensoring with the rationals,
$\pi_* (G) \otimes \QQ$ with the induced
bracket is a graded Lie algebra over $\QQ.$

The following result can be deduced in various ways;
our proof is based on results of H.~Scheerer in~\cite{Sch}.

\begin{Thm}\label{abelian}
Let $(G, \mu)$ be a connected CW~homotopy-associative H-space.
The following are equivalent.
\begin{enumerate}
\item\label{G1}
$(G, \mu)$ has the rational H-type of an abelian topological group.
\item\label{G2}
$(G, \mu)$ is rationally homotopy-abelian.
\item\label{G3}
The Samelson Lie algebra $\pi_{*} (G) \otimes \QQ,$
described above, is abelian, that is, all brackets are zero.
\item\label{G4}
The rational equivalence
\[
e \colon G \to \prod_{j \geq 1} K (\pi_{j} (G) \otimes \QQ, j)
\]
of Lemma~\ref{G decomposed} is a rational H-equivalence,
where the product of Eilenberg-Mac Lane
spaces has the standard multiplication.
\end{enumerate}
\end{Thm}

\begin{proof}
The implication $(\ref{G1}) \!\implies \! (\ref{G2})$
is immediate from definitions.
The implication $(\ref{G4}) \!\implies \! (\ref{G1})$
is the result of Milnor mentioned in Example~\ref{EMspace}
(Corollary to Theorem~3 of~\cite{Mil2}).

We prove $(\ref{G2}) \!\implies \! (\ref{G3}).$
Proposition~\ref{James} implies that $(G, \mu)$ is group-like,
so $(G_{\QQ}, \mu_{\QQ})$ is group-like too.
Corollary~2 in Section~0.1 of~\cite{Sch} therefore shows that
$(G, \mu)$ is rationally H-equivalent to a loop space $\Omega X$
with the usual multiplication of loops.
Now use Theorem~\ref{nil<=Hnil}.

Finally, for $(\ref{G3}) \!\implies \! (\ref{G4})$
use Corollary~3 of~\cite{Sch},
which gives a bijection between the set of H-maps
between two rational group-like spaces
and the space of homomorphisms between their rational Samelson
algebras.
The rational Samelson algebra
of $\prod_{j \geq 1} K (\pi_{j} (G), j)$ is abelian because
$\prod_{j \geq 1} K (\pi_{j} (G), j)$ is rationally abelian.
So if $\pi_{*} (G) \otimes \QQ$ is abelian then $e$ trivially induces
an isomorphism of rational Samelson algebras and is thus a
rational H-equivalence.
\end{proof}

The following useful consequences are well-known.

\begin{Cor}\label{P_Loop}
The standard loop multiplication on a product of Eilenberg Mac Lane
spaces is the unique group-like, homotopy-abelian H-structure
up to rational H-equivalence.
\end{Cor}

\begin{proof}
Use (\ref{G2}) $\iff$ (\ref{G4}) in Theorem~\ref{abelian}.
\end{proof}

\begin{Cor} \label{G abelian}
Let $G$ be a connected CW~homotopy-associative H-space.
Suppose $H^* (G; \QQ)$ is finite-dimensional.
Then $G$ is rationally homotopy-abelian.
\end{Cor}

\begin{proof}
Referring to the proof of Lemma~\ref{G decomposed},
we see that finite-dimensionality of $H^* (G; \QQ)$
implies that $H^* (G; \QQ)$ is an exterior algebra on
odd-degree generators.
By that same result, we conclude that $\pi_{*} (G) \otimes \QQ$
is zero in even degrees.
Thus, for degree reasons,
the rational Samelson bracket on $\pi_{*} (G) \otimes \QQ$ is trivial.
This is Condition~(\ref{G3}) of Theorem~\ref{abelian},
so we conclude that $G$ is rationally homotopy-abelian.
\end{proof}

Finally, we identify the H-spaces $F (X, G)\idc$
and $F_{\bullet} (X, G)\idc,$
for $X$ and $G$ finite complexes, up to rational H-equivalence.
Note that a continuous map $f \colon X \to Y$
induces H-maps $f^{*} \colon F (Y, G)\idc \to F (X, G)\idc$
and $f^{*} \colon F_{\bullet} (Y, G)\idc \to F_{\bullet} (X, G)\idc.$
Similarly, an H-map $h \colon G \to H$ induces H-maps
$h_{*} \colon F (X, G)\idc \to F (X, H)\idc$ and
$h_{*} \colon F_{\bullet} (X, G)\idc \to F_{\bullet} (X, H)\idc.$

\begin{Thm}\label{F(X,G)}
Let $X$ be a finite CW~complex.
Let $G$ be a connected CW~homotopy-associative H-space
with $H^* (G; \QQ)$ finite-dimensional.
Then $F (X, G)\idc$ and $F_{\bullet} (X, G)\idc$ are
homotopy-abelian after rationalization.
Consequently, the rational equivalences of Theorem~\ref{null-maps}
are actually rational H-equivalences
where the products of Eilenberg-Mac Lane spaces have
the standard multiplication.
Moreover, these rational H-equivalences are
natural with respect to maps $f \colon X \to Y$
and H-maps $h \colon G \to H.$
\end{Thm}

\begin{proof}
Corollary~\ref{G abelian} implies that
$G$ is rationally homotopy-abelian.
Thus, by Theorem~\ref{Hnil},
$F (X, G_{\QQ})\idc$ is rationally homotopy-abelian.
Using Proposition~\ref{Samequiv5} and the fact that the rational
Samelson algebra of an H-product is the product of the
Samelson algebras (see Example~21.4 of~\cite{FHT}), we see
that the rational Samelson algebra of $F_{\bullet} (X, G)\idc$ is
also abelian.
Thus $F_{\bullet} (X, G)\idc$ is
rationally homotopy-abelian by Theorem~\ref{abelian}.
The naturality assertions now follow from the naturality
of the rational homotopy equivalence of a rational H-space
with the appropriate product of Eilenberg-Mac Lane spaces as in
Corollary~\ref{G_Dec_Cor}.
\end{proof}

\section{Conclusion: Passage to Limits}\label{sec:conclusion}

We now prove the main results of the paper.
First, using our results above we give a preliminary version
of Theorem~\ref{intromain1}
in the special case in which
the maximal ideal space $\MM (A)$ of $A$
happens to be a finite complex.

\begin{Thm}\label{rational-finitecomplex}
Suppose that $A$ is a commutative unital Banach algebra
and $\MM (A)$ is a finite CW~complex.
Let $V_{h, j} = H^{2 j - 1 - h} (\MM (A); \QQ).$
Then there is a rational H-equivalence
\[
\GL_n (A)\idc \, \simeq_{\QQ} \,
\prod_{j = 1}^n
\prod_{h = 1}^{2 j - 1} K (V_{h, j}, h),
\]
\item
where the product of Eilenberg-Mac Lane
spaces has the standard loop multiplication.
The equivalence is natural with respect to homomorphisms
between commutative
unital Banach algebras with maximal ideal space a finite complex.
Moreover, for fixed $A$ with $\MM (A)$ a finite complex
the equivalence is natural in $n$ in the following sense:
For each $n > 1,$ let
$i \colon \GL_{n - 1} (A)\idc \to \GL_n (A)\idc$ denote the
inclusion in the upper left corner.
Then, after the identifications
\[
\pi_k (\GL_{n - 1} (A)\idc) \otimes \QQ
  \cong \bigoplus_{j = 1}^{n - 1} V_{k, j}
\hbox {\ \ and \ \ }
\pi_k (\GL_{n} (A)\idc) \otimes \QQ
 \cong  \bigoplus_{j =1 }^{n} V_{k, j}
\]
for each $k \geq 1$ given above,
the map
\[
i_{\sharp} \otimes 1 \colon \pi_k (\GL_{n - 1} (A)\idc) \otimes \QQ
  \to \pi_k (\GL_{n} (A)\idc) \otimes \QQ
\]
corresponds to the inclusion of vector spaces
\[
\bigoplus_{j = 1}^{n - 1} V_{k, j}
  \hookrightarrow \bigoplus_{j = 1}^{n} V_{k, j}.
\]
\end{Thm}

\begin{proof}
By Theorem~\ref{Davie1} we have an H-equivalence
\[
\GL_n (A)\idc \simeq F (\MM (A), \GL_n (\CC))\idc.
\]
The rational homotopy groups of $\GL_n (\CC) \simeq U_n (\CC)$
occur in degrees $1, 3, \dots, 2 n - 1$ by Example~\ref{U_n (C)}.
Theorem~\ref{null-maps} thus gives the needed rational equivalence.
Since $ H^{*} (\GL_{n} (\CC), \QQ)$ is finite-dimensional,
Theorem~\ref{F(X,G)} implies this equivalence is actually a
rational H-equivalence.
Naturality with respect to Banach algebra homomorphisms
is a direct consequence
of the naturality given in that theorem.
For naturality with respect to $n,$
we first need to know what the inclusion
$\GL_{n - 1} (\CC)\idc \hookrightarrow \GL_n (\CC)\idc$
does on the rational homotopy groups.
For this, use the retractions to the corresponding unitary
groups and Example~\ref{U_n (C)}.
Naturality with respect to $n$ now follows from the naturality in
Theorem~\ref{F(X,G)}
and the fact that the inclusion
$i \colon \GL_{n - 1} (A)\idc \to \GL_n (A)\idc$
corresponds, via the H-equivalence of Theorem~\ref{Davie1},
to the map
$F (\MM (A), \GL_{n - 1} (\CC)) \to F (\MM (A), \GL_{n} (\CC))\idc$
induced by the inclusion
$\GL_{n - 1} (\CC)\idc \hookrightarrow \GL_n (\CC)\idc.$
\end{proof}

We deduce the general case by considering limits of direct
systems in our various settings.
Recall that a unital homomorphism $\ph \colon A \to B$
between unital Banach algebras
induces a group homomorphism
$\GL_n (\ph) \colon \GL_n (A) \to \GL_n (B).$
We first consider various direct systems in the \ca\   setting.
The material here is in principle well known,
but we have been unable to find a reference giving it in
the generality needed here.
The closest we are aware of is Section~4 of~\cite{Hn}.
Unfortunately, in this reference it is assumed
throughout that all the maps in the direct system are injective,
a condition which does not hold in our application.

Let $\Ld$ be a directed set.
Let $(A_{\ld})_{\ld \in \Ld}$ be a direct system of \ca s,
indexed by $\Ld$
with \ca\  maps $\ph_{\ld, \mu} \colon A_{\ld} \to A_{\mu}.$
The direct limit $A = \Dirlim_{\ld} A_{\ld}$
exists in the category of \ca s.
In fact, we have the following characterization.

\begin{Thm} \label{Cdirlim}
The direct limit $A = \Dirlim_{\ld} A_{\ld}$ of a system
of \ca\  maps $\ph_{\ld, \mu} \colon A_{\ld} \to A_{\mu}$
is characterized by the existence of \ca\  homomorphisms
$\ph_{\ld, \infty} \colon A_{\ld} \to A$ such that:
\begin{enumerate}
\item\label{P_DLimDfnCons}
$\ph_{\mu, \infty} \circ \ph_{\ld, \mu} = \ph_{\ld, \infty}$
whenever $\mu \geq \ld.$
\item\label{P_DLimDfnNorm}
If $a \in A_{\ld}$ and $\ph_{\ld, \infty} (a) = 0$ then
$\lim_{\mu} \| \ph_{\ld, \mu} (a) \| = 0.$
\item\label{P_DLimDfnDense}
$\bigcup_{\ld \in \Ld} \ph_{\ld, \infty} (A_{\ld})$ is dense in $A.$
\end{enumerate}
\end{Thm}

\begin{proof}
See the proof of Proposition 2.5.1 of~\cite{Ph1}
for the construction,
and see Theorem 6.1.2 of~\cite{Mr} and the preceding discussion
for the case $\Ld = {\mathbb{N}}.$
That these properties characterize the direct limit
for $\Ld = {\mathbb{N}}$
is implicit in the proof of Theorem 6.1.2(b) of~\cite{Mr},
and the proof for a general index set is the same.
(Murphy has a slightly different formulation of the second condition.
The key fact in the proof is that
an injective *-algebra homomorphism of \ca s is
isometric, even without assuming continuity to begin with.)
\end{proof}

We emphasize that $\bigcup_{\ld \in \Ld} \ph_{\ld, \infty} (A_{\ld})$
is usually {\emph{not}} the entire direct limit,
merely a dense subalgebra.

Given a direct system $(A_{\ld})_{\ld \in \Ld}$ of \ca s
with maps $\ph_{\ld, \mu} \colon A_{\ld} \to A_{\mu},$
we obtain direct systems of matrix and function spaces.
Recall $M_n (A)$ denotes the space of $n \times n$
matrices with entries in $A.$
Write
$\widetilde{\ph}_{\lambda, \mu}
 \colon M_n (A_{\lambda}) \to M_n (A_{\mu})$
for the map induced by applying $\ph_{\lambda, \mu}$
entry-wise.
For function spaces in the \ca\  context,
we introduce the following notation.
We deviate from our usual notation
for function spaces to follow the standard notation
$C (Y, A)$ for \ca s.

\begin{Ntn} \label{C (Y, B)}
For any \ca\  $A$ and any space $Y$
let $C (Y, A)$ denote the \ca\  of continuous functions
$b \colon Y \to A,$ with pointwise operations and the
supremum norm $\| a \| = \sup_{y \in Y} \| a (y) \|.$
If moreover $\ph \colon A \to B$ is a homomorphism of \ca s,
we write
$\overline{\ph} \colon C (Y, A) \to C (Y, B)$
for the map obtained by composition with $\ph.$
\end{Ntn}

The following result shows the direct limits in both cases
are as expected.
The result is a special case of the statement
that direct limits commute with maximal completed tensor products,
but our proof avoids the
technicalities of completed tensor products.

\begin{Lem}\label{P_DLimTens}
Let $(A_{\ld})_{\ld \in \Ld}$ be a direct system of \ca s,
with maps $\ph_{\ld, \mu} \colon A_{\ld} \to A_{\mu},$
and let $A = \Dirlim_{\ld} A_{\ld}.$
Let $Y$ be a compact Hausdorff space.
Then:
\begin{enumerate}
\item\label{P_DLimTensMn}
For $n \in {\mathbb{N}}$ there is a natural isomorphism
$\Dirlim_{\ld} M_n (A_{\ld}) \cong M_n (A).$
\item\label{P_DLimTensCX}
Then there is a natural isomorphism
$\Dirlim_{\ld} C (Y, A_{\ld}) \cong C (Y, A).$
\end{enumerate}
\end{Lem}

\begin{proof}
For~(\ref{P_DLimTensMn}),
the maps are obtained by applying $\ph_{\ld, \infty}$ to each entry.
In the characterization of the direct limit (Theorem~\ref{Cdirlim}),
Parts~(\ref{P_DLimDfnCons}) and~(\ref{P_DLimDfnDense})
are immediate, and Part~(\ref{P_DLimDfnNorm}) follows from
the fact that if the entries of a net of matrices converge
to zero, then so do the matrices.

We prove~(\ref{P_DLimTensCX}).
The maps $C (Y, A_{\ld}) \to C (Y, A_{\mu})$
and $C (Y, A_{\ld}) \to C (Y, A)$
are of course ${\overline{\ph}}_{\ld, \mu}$
and ${\overline{\ph}}_{\ld, \infty},$
as in Notation~\ref{C (Y, B)}.

Part~(\ref{P_DLimDfnCons})
of the characterization of the direct limit is immediate.
For Part~(\ref{P_DLimDfnNorm}),
let $b \in C (Y, A_{\ld})$
with ${\overline{\ph}}_{\ld, \infty} (b) (y) = 0.$
Because homomorphisms of \ca s are norm decreasing,
the functions $f_{\mu} (y) = \| \ph_{\ld, \mu} (b (y)) \|$
decrease pointwise to~$0.$
By Dini's Theorem (Problem~E in Chapter~7 of~\cite{Kl}),
the convergence is uniform,
that is, $\lim_{\mu} \| {\overline{\ph}}_{\ld, \mu} (b) \| = 0.$

We finish by proving Part~(\ref{P_DLimDfnDense})
of the characterization of the direct limit.
Let $b \in C (Y, A),$ and let $\ep > 0.$
Choose a finite partition of unity $(f_1, \ldots, f_n)$
such that whenever $x, y \in {\mathrm{supp}} (f_k)$
then $\| b (x) - b (y) \| < \frac{1}{3} \ep.$
Choose $y_k \in {\mathrm{supp}} (f_k).$
Choose $\ld_k \in \Ld$ and $c_k \in A_{\ld_k}$ such that
$\| \ph_{\ld_k, \infty} (c_k) - b (y_k) \| < \frac{1}{3} \ep.$
Choose $\ld \in \Ld$ such that $\ld \geq \ld_k$ for all $k,$
and let $a_k = \ph_{\ld_k, \ld} (c_k) \in A_{\ld}.$
Then $\| \ph_{\ld, \infty} (a_k) - b (y_k) \| < \frac{1}{3} \ep.$
Define $a \in A_{\ld}$ by $a (y) = \sum_{k = 1}^n f_k (y) a_k.$
For $y \in Y$ we have $\| b (y_k) - b (y) \| < \frac{1}{3} \ep$
whenever $f_k (y) \neq 0,$ whence
\begin{align*}
\| \ph_{\ld, \infty} (a (y)) - b (y) \|
& \leq \sum_{k = 1}^n f_k (y)
      \left[ \rule{0em}{2.0ex}
         \| a_k - b (y_k) \| + \| b (y_k) - b (y) \| \right]  \\
& \leq \sum_{k = 1}^n 2 \left( \frac{\ep}{3} \right) f_k (y)
  = \frac{2 \ep}{3}.
\end{align*}
So $\| {\overline{\ph}}_{\ld, \infty} (a) - b \| < \ep.$
\end{proof}

Now given a direct system $(A_{\ld})_{\ld \in \Ld}$ of unital \ca s
with unital maps $\ph_{\ld, \mu} \colon A_{\ld} \to A_{\mu},$
we obtain a direct system $(\GL_n (A_{\ld}))_{\ld \in \Lambda}$
of invertible groups
with structure maps
$\GL_n (\ph_{\ld, \mu}) \colon \GL_n (A_{\ld}) \to \GL_n (A_{\mu})$
The following result is implicit in the usual direct proofs that
$K_1$ commutes with direct limits of \ca s.
The proofs we know in the literature, however,
instead use Bott periodicity and the result for $K_0.$

\begin{Thm}\label{P_DLimPi_n}
Let $(A_{\ld})_{\ld \in \Ld}$ be a direct system of unital \ca s,
with unital maps $\ph_{\ld, \mu} \colon A_{\ld} \to A_{\mu},$
and let $A = \Dirlim_{\ld} A_{\ld}.$
Then the maps $\GL_n (\ph_{\ld, \infty}) \colon \GL_n (A_{\ld}) \to
\GL_n (A)$ induce a natural isomorphism
$\Dirlim_{\ld} \pi_k (\GL_n (A_{\ld}) ) \cong \pi_k (\GL_n (A))$
for $k \geq 0.$
\end{Thm}

\begin{proof}
Following Notation~\ref{C (Y, B)},
set $B_{\ld} = C (S^k, M_n (A_{\ld}))$ and $B = C (S^k, M_n (A)).$
Then $\pi_k (\GL_n (A_{\ld}) ) = \pi_0 (\GL_1 (B_{\ld}))$
and $\pi_k (\GL_n (A) ) = \pi_0 (\GL_1 (B)),$
and Lemma~\ref{P_DLimTens}
identifies $B$ naturally with $\Dirlim_{\ld} B_{\ld}.$
Therefore it suffices to prove the lemma when $k = 0$ and $n = 1.$

We must prove:
\begin{enumerate}
\item\label{P_DLimPi_n_Surj}
Every path component of $\GL_1 (A)$ contains some element
$\ph_{\ld, \infty} (b)$
for some $\ld$ and some $b \in \GL_1 (A_{\ld}).$
\item\label{P_DLimPi_n_Inj}
If $b \in \GL_1 (A_{\ld})$ and
$\ph_{\ld, \infty} (b) \in \GL_1 (A)\idc,$
then there exists $\mu \geq \ld$
such that $\ph_{\ld, \mu} (b) \in \GL_1 (A_{\mu})\idc.$
\end{enumerate}

The key ingredient is that if $s \in A$ is invertible,
then so is any $t \in A$ with $\| t - s \| < \| s^{-1} \|^{-1}.$
We use this to prove the following claim:
if $b \in A_{\ld}$ and
$\ph_{\ld, \infty} (b)$ is invertible, then there exists
$\mu \geq \ld$ such that $\ph_{\ld, \mu} (b)$ is invertible.
Let $a = \ph_{\ld, \infty} (b).$
Set $\ep_0 = 1 / (3 \| a \|).$
Choose $\ld_0 \in \Ld$ and $c_0 \in A_{\ld_0}$ such that
$\| \ph_{\ld_0, \infty} (c_0) - a^{-1} \| < \ep_0.$
Then
\[
\| a \ph_{\ld_0, \infty} (c_0) - 1 \|
  < \ep_0 \| a \|
  \leq \tfrac{1}{3}.
\]
Choose $\ld_1 \geq \ld, \ld_0.$
By
Part~(\ref{P_DLimDfnNorm}) of the characterization of the direct limit,
there exists $\mu_0 \geq \ld_0$ such that
\[
\| \ph_{\ld, \mu_0} (b) \ph_{\ld_0, \mu_0} (c_0) - 1 \| < \tfrac{2}{3}.
\]
In particular,
$\ph_{\ld, \mu_0} (b) \ph_{\ld_0, \mu_0} (c_0)$ is invertible.
Similarly, there is $\mu \geq \mu_0$ such that
$\ph_{\ld_0, \mu} (c_0) \ph_{\ld, \mu} (b)$ is invertible.
Then $\ph_{\ld_0, \mu} (b)$ is invertible.
The claim is proved.

To prove~(\ref{P_DLimPi_n_Surj}), let $a \in \GL_1 (A),$
and use density of
$\bigcup_{\ld \in \Ld} \ph_{\ld, \infty} (A_{\ld})$ in $A$
to choose $\ld_0$ and $b \in A_{\ld_0}$ such that
$\| \ph_{\ld_0, \infty} (b) - a \| < \| a^{-1} \|^{-1}.$
Then the straight line path
from $\ph_{\ld_0, \infty} (b_0)$ to $a$ is in $\GL_1 (A).$
By the claim, there is $\ld \geq \ld_0$ such that
$b = \ph_{\ld_0, \ld} (b_0)$ is invertible.

To prove~(\ref{P_DLimPi_n_Inj}), let $b \in \GL_1 (A_{\ld})$
satisfy $\ph_{\ld, \infty} (b) \in \GL_1 (A)\idc.$
Let $t \mapsto a_0 (t)$
be a continuous path from $\ph_{\ld, \infty} (b)$
to~$1$ in $\GL_1 (A),$
defined for $t \in [ 1, 2].$
Regard $t \mapsto a_0 (t)$ as an invertible element
$a \in C ( [ 1, 2], A ).$
By Lemma~\ref{P_DLimTens}(\ref{P_DLimTensCX}),
we have
$C ( [ 1, 2], A ) = \Dirlim_{\ld} C ( [ 1, 2], A_{\ld} ),$
so there is $\ld_0 \in \Ld$ and $c_0 \in C ( [ 1, 2], A_{\ld_0})$
such that
$\| {\overline{\ph}}_{\ld_0, \infty} (c_0) - a_0 \|
                    < \| a_0^{-1} \|^{-1}.$
Without loss of generality $\ld_0 \geq \ld.$
Let $a \in C ([0, 3], A)$ be the concatenation
of the constant path $\ph_{\ld, \infty} (b)$ on $[ 0, 1],$
the path $a_0$ on $[ 1, 2],$
and the constant path~$1$ on $[ 2, 3].$
Let $c \in C ([0, 1], A_{\ld_0})$ be the concatenation
of the straight line path from $\ph_{\ld, \ld_0} (b)$
to $c_0 ( 1)$ on $[ 0, 1],$
the path $c_0$ on $[ 1, 2],$
and the straight line path from $c_0 ( 2)$ to~$1$ on $[ 2, 3].$
Then also
$\| {\overline{\ph}}_{\ld_0, \infty} (c) - a \| < \| a^{-1} \|^{-1}.$
So ${\overline{\ph}}_{\ld_0, \infty} (c)$ is invertible.
Applying the claim to the direct system
$(C ([0, 3], A_{\ld}))_{\ld \in \Ld}$
(using Lemma~\ref{P_DLimTens}(\ref{P_DLimTensCX}) again),
we find $\mu \geq \ld_0$ such that
${\overline{\ph}}_{\ld_0, \mu} (c)$ is invertible.
Since $\ph_{\ld_0, \mu} (c (3)) = 1,$ we get
$\ph_{\ld, \mu} (b) = \ph_{\ld_0, \mu} (c (0)) \in \GL_1 (A_{\mu})\idc.$
\end{proof}

Finally, we briefly consider direct limits
of graded Lie algebras over $\QQ.$
Let $\Ld$ be a directed set,
for $\ld \in \Ld$ let $L_{\ld}$ be a graded Lie
algebra over $\QQ$ as in Definition~\ref{gla}, and suppose given
graded Lie algebra maps $\ps_{\ld, \mu} \colon L_{\ld} \to L_{\mu}$
for $\ld \leq \mu$ satisfying the usual coherence conditions.
Define the direct limit $\Dirlim_{\ld} L_{\ld}$ to be
the graded space given, in each degree $n > 0,$
as the algebraic direct limit of the
$L_{\ld}$ in degree $n$ with the bracket induced by the maps
$\ps_{\ld, \infty} \colon L_{\ld} \to \Dirlim_{\ld} L_{\ld}.$
It is direct to check $L$ with the induced bracket satisfies the
axioms in Definition~\ref{gla}.
We need one easy result in this context.

\begin{Thm}\label{lim-gla}
Let $(L_{\ld})_{\ld \in \Lambda}$ be a direct system of graded Lie
algebras over $\QQ$ with graded Lie algebra structure maps
$\ps_{\ld, \mu} \colon L_{\ld} \to L_{\mu}.$
The direct limit $\Dirlim_{\ld} L_{\ld}$ satisfies:
\begin{enumerate}
\item
Given a graded Lie algebra $L$ over $\QQ$
and graded Lie algebra
maps $\ph_{\ld} \colon L_{\ld} \to L$ satisfying
$\ph_{\mu} \circ \ps_{\ld, \mu} = \ph_{\ld}$
there exists a unique graded Lie algebra map
$\ph \colon \Dirlim_{\ld} L_{\ld} \to L$
satisfying $\ph \circ \ps_{\ld, \infty} = \ph_{\ld}.$
\item
If each $L_{\ld}$ is abelian then so is $\Dirlim_{\ld} L_{\ld}.$
\end{enumerate}
\end{Thm}

\begin{proof}
The proof of~(1) is a direct consequence of the result for direct
limits of groups.
The proof of~(2) is direct from the definition of the bracket
in $\Dirlim_{\ld} L_{\ld}.$
\end{proof}

We can now give the proof of our first main result.

\begin{proof}[Proof of Theorem~\ref{intromain1}]
By Theorem~\ref{Davie1},
it suffices to prove the rational H-equivalence for
algebras of the form $C (X),$ with $X$ compact Hausdorff,
and to prove naturality for unital homomorphisms between
such algebras,
that is, for continuous maps (in the opposite direction) of
compact Hausdorff spaces.

So let $X$ be compact Hausdorff.
By Theorem~10.1 in Chapter~X of~\cite{ES},
there exists a directed
set $\Ld$ and an inverse system
$(X_{\ld})_{\ld \in \Ld}$ of finite CW~complexes and continuous maps
$f_{\ld, \mu} \colon X_{\mu} \to X_{\ld}$
for $\ld \leq \mu$ with
$X \approx \Invlim_{\ld} X_{\ld}.$
We then obtain a direct system of \ca s
$(C (X_{\ld}))_{\ld \in \Lambda}$
with maps $\ph_{\ld, \mu} \colon C (X_{\ld}) \to C (X_{\mu})$
induced by the $f_{\ld, \mu}.$
The maps $X \to X_{\ld}$
determine an isomorphism of \ca s $\Dirlim_{\ld} C (X_{\ld}) \to C (X)$
because $Y \mapsto C (Y)$ is a contravariant
equivalence of categories from compact Hausdorff spaces
to commutative unital \ca s.
Applying Theorem~\ref{P_DLimPi_n}
and restricting to the identity components,
we obtain isomorphisms
$\Dirlim_{\ld} \pi_k (\GL_n (C (X_{\ld}))\idc) \cong
      \pi_k (\GL_n (C (X))\idc)$
for $k > 0.$
Since the tensor product commutes with direct limits, we obtain
\[
\Dirlim_{\ld} \pi_k (\GL_n (C (X_{\ld}))\idc) \otimes \QQ
 \cong \pi_k (\GL_n (C (X))\idc) \otimes \QQ.
\]

Set $V^{\ld}_{k, j} = H^{2 j - 1 - k} (X_{\ld}; \QQ)$
and let
\[
H (f_{\ld, \mu}) \colon \bigoplus_{j = 1}^{n} V^{\ld}_{k, j} \to
\bigoplus_{j = 1}^{n}V^{\mu}_{k, j}
\]
be the map induced on rational cohomology groups by $f_{\ld, \mu}.$
Since each $\MM (C (X_{\ld})) \approx X_{\ld}$
is a finite complex,
Theorem~\ref{rational-finitecomplex} applies to
give a commutative diagram
\[
\xymatrix{ \pi_k (\GL_n (C (X_{\ld}))\idc) \otimes \QQ
\ar[rr]^{\ \ \ \ \ \ \cong}
\ar[d]_{\GL_n (\ph_{\ld, \mu})_{\sharp} \otimes 1}
&& \bigoplus_{j = 1}^{n} V^{\ld}_{k, j}
\ar[d]^{H (f_{\ld, \mu})}\\
\pi_k (\GL_n (C (X_{\mu}))\idc) \otimes \QQ
\ar[rr]^{\ \ \ \ \ \ \cong}
&& \bigoplus_{j = 1}^{n} V^{\mu}_{k, j}}
\]
for each $k > 0.$
By the
continuity of \v{C}ech cohomology
(Theorem~3.1 in Chapter~X of~\cite{ES};
see page~110 for the definition of ${\mathcal{G}}_R$),
$\Dirlim_{\ld} \, V^{\ld}_{k, j}
 \cong {\Check{H}}^{2 j - 1 - k} (X; \QQ) = {\Check{V}}_{k, j}.$
We conclude
\[
\pi_k (\GL_n (C (X))\idc) \otimes \QQ \cong
\bigoplus_{j = 1}^{n} {\Check{V}}_{k, j}.
\]

Naturality is clear.

It remains to prove the rational H-equivalence and show that
$\GL_n (C (X))\idc$ is rationally homotopy-abelian.
Consider the direct system
$(\pi_* (\GL_n (C (X_{\ld}))\idc) \otimes \QQ)_{\ld \in \Ld}$
of graded Lie algebras over $\QQ$ with structure maps
\[
\GL_n (\ph_{\ld, \mu}) \otimes 1
  \colon \pi_* (\GL_n (C (X_{\ld}))\idc) \otimes \QQ
  \to \pi_* (\GL_n (C (X_{\mu}))\idc) \otimes \QQ.
\]
By Theorems~\ref{P_DLimPi_n} and~\ref{lim-gla}(1) there is
a graded Lie algebra isomorphism
\[
\Dirlim_{\ld}
\pi_* (\GL_n (C (X_{\ld}))\idc)
        \otimes \QQ \cong \pi_* (\GL_n (C (X))\idc) \otimes \QQ.
\]
By Theorems~\ref{rational-finitecomplex}
and~\ref{abelian}(\ref{G3}),
each $\pi_* (\GL_n (C (X_{\ld}))\idc) \otimes \QQ$ is abelian.
Therefore $\pi_* (\GL_n (C (X))\idc)\otimes \QQ$ is abelian by
Theorem~\ref{lim-gla}(2).
The space $\GL_n (C (X))$ is homotopy equivalent to a CW~complex
(Theorem~\ref{CR_Serre}(\ref{P_CR5})).
So we may apply
$(\ref{G3}) \!\implies \! (\ref{G2})$
and $(\ref{G3}) \!\implies \! (\ref{G4})$
of Theorem~\ref{abelian}.
\end{proof}

\begin{Cor} \label{GL-injection}
Let $A$ be a commutative unital Banach algebra and $n > 1.$
Let
$i \colon \GL_{n - 1} (A)\idc \to \GL_n (A)\idc$ denote the
inclusion.
Let ${\Check{V}}_{h, j} = {\Check{H}}^{2 j - 1 - h} (\MM (A); \QQ).$
Then, after the identifications
\[
\pi_k (\GL_{n - 1} (A)\idc) \otimes \QQ
 \cong \bigoplus_{j = 1}^{n - 1} {\Check{V}}_{k, j}
\hbox {\ \ and \ \ } \pi_k (\GL_{n} (A)\idc) \otimes \QQ
 \cong \bigoplus_{j = 1}^{n} {\Check{V}}_{k, j}
\]
for each $k \geq 1$ given by Theorem~\ref{intromain1},
the map
\[
i_{\sharp} \otimes 1 \colon \pi_k (\GL_{n - 1} (A)\idc) \otimes
\QQ \to \pi_k (\GL_{n} (A)\idc) \otimes \QQ
\]
corresponds to the inclusion of vector spaces
\[
\bigoplus_{j = 1}^{n - 1} {\Check{V}}_{k, j} \hookrightarrow
\bigoplus_{j = 1}^{n} {\Check{V}}_{k, j}.
\]
\end{Cor}

\begin{proof}
As in the proof of Theorem~\ref{intromain1} above,
we may assume $A = C (X)$ and write
\[
\pi_k (\GL_{n - 1} (A)\idc) \otimes \QQ
 = \Dirlim_{\ld}\pi_k (\GL_{n - 1} (C (X_{\ld}))\idc)
  \otimes \QQ
\]
and
\[
\pi_k (\GL_{n} (A)\idc) \otimes \QQ
 = \Dirlim_{\ld}\pi_k (\GL_{n - 1} (C (X_{\ld}))\idc)
  \otimes \QQ
\]
where each $X_{\ld}$ is a finite complex.
Now use naturality with respect to $n$ in
Theorem~\ref{rational-finitecomplex}.
\end{proof}

We can now deduce Theorem~\ref{intromain2}.

\begin{proof}[Proof of Theorem~\ref{intromain2}]
Consider the long exact sequence of Theorem~\ref{T:les1}
as far as the last abelian term $\pi_1 (\GL_n (A)).$
In the groups $\pi_k (-)$ for $k \geq 1,$
we may replace each space by the component containing the basepoint,
obtaining
\[
\xymatrix{
\cdots \ar[r] & \pi_k ( \GL_{n - 1} (A)\idc ) \ar[r]
 & \pi_k ( \GL_{n} (A)\idc ) \ar[r] & \pi_k ( \Lc_n (A)\idc)
  \ar `d[lll] `[llld]_{}[llld] \\
\pi_{k - 1} (\GL_{n - 1} (A)\idc) \ar[r] &
  \ \ \ \   \cdots  \ \ \ \ \ar[r] &
 \pi_1 (\GL_n (A)\idc) \ar[r] & \pi_1 (\Lc_n (A)\idc).}
\]
Tensor with $\QQ.$
Applying Corollary~\ref{GL-injection}, we compute
\[
\pi_k (\Lc_n (A)\idc) \otimes \QQ
 \cong {\Check{H}}^{2 n - 1 - k} (\MM (A); \QQ)
\]
for $k > 1,$ as needed.
\end{proof}

We conclude with some remarks and a result
concerning the rational homotopy type of
the space of last columns for a commutative unital Banach algebra.
By Theorem~\ref{commutativeLcn}, we have $\Lc_n (A)\idc \simeq
F (\MM (A), S^{2 n - 1})\idc.$
Thus $\Lc_n (A)\idc$ is a nilpotent CW~complex
when $\MM (A)$ is a finite complex by Theorem~\ref{HMRlocal}.
Using Example~\ref{odd-sphere} and
Theorem~\ref{null-maps},
we conclude in this case that $\Lc_n (A)_{\circ}$
is a rational H-space.
When $\MM (A)$ is not a finite complex, there is no guarantee
that $\Lc_n (A)\idc$ is a nilpotent space
(although we know of no counterexamples),
and we cannot discuss its rational homotopy type.
However, it follows from Theorem~\ref{CR_Serre}(\ref{P_CR5})
that $\Lc_n (A)\idc$
has the homotopy type of a CW~complex and thus admits a universal
cover.
Our last result describes the rational homotopy type of the
universal cover of $\Lc_n (A)\idc$ in the general case of a
commutative unital Banach algebra.
We need one final result from rational homotopy theory.

\begin{Lem}\label{ratsur}
Let $p \colon E \to B$ be a map of nilpotent CW~complexes inducing a
surjection on rational homotopy groups.
If $E$ is a rational H-space
then $B$ is one also.
\end{Lem}

\begin{proof}
This result is well-known in rational homotopy theory and easy to
prove using minimal models.
Since we have avoided their use thus
far, we give a proof which does not use minimal models.
Our argument adapts one
given in~\cite{Sch}; we reproduce it here for the sake of
completeness.

As usual,
we assume that the basepoint of each CW~H-space $(G, \mu)$
is the identity.
In addition, following Theorem III.4.7 of~\cite{Wh},
we assume that the identity $e_G$ is strict,
in the sense that $\mu (x, e_G) = \mu (e_G, x) = x$
for all $x \in G.$

Write $\pi_* (B_{\QQ}) = \bigoplus_{i \in J} V_i$ with each $V_i$ a
$1$-dimensional rational vector space concentrated in degree $n_i.$
For each $i,$ pick a basis element $\alpha_i$ of $V_i.$
Since $(p_{\QQ})_{\#} \colon \pi_* (E_{\QQ}) \to \pi_* (B_{\QQ})$
is surjective, we
may choose $\beta_i \in \pi_{n_i} (E_{\QQ})$ such that
$(p_{\QQ})_{\#} (\beta_i) = \alpha_i$ for each $i.$
If $n_i$ is odd,
then $K (\QQ, n_i) \simeq S^{n_i}_{\QQ}$
and so we may regard $\beta_i$ as
a map $\gamma_i \colon K (\QQ, n_i) \to E_{\QQ}.$
Now suppose that $n_i$ is even.
Then $K (\QQ, n_i) \simeq \Omega \Sigma S^{n_i}_{\QQ}.$
For any space $X,$
let $\ep_X \colon X \to \Omega \Sigma X$
be the adjoint of the suspension of the identity map of $X.$
This map is natural in $X.$
Since $E_{\QQ}$ is an H-space, by Theorem~8.14 of~\cite{St}
there is a retraction
$r \colon \Omega \Sigma E_{\QQ} \to E_{\QQ}$
of $\ep_{E_{\QQ}}.$
Using this retraction, we have a commutative diagram
\[
\xymatrix{ K (\QQ, n_i) \ar[r]^{\simeq}
 & \Omega \Sigma S^{n_i}_{\QQ} \ar[r]^-{\Omega \Sigma \beta_i}
 & \Omega \Sigma E_{\QQ} \ar[d]^{r} \\
 & S^{n_i}_{\QQ} \ar[u]^{\varepsilon_{S^{n_i}_{\QQ}}} \ar[r]_-{\beta_i}
 & E_{\QQ}.}
\]
Now set
$\gamma_i
 = r \circ \Omega \Sigma \beta_i \colon K (\QQ, n_i) \to E_{\QQ}.$

Choose a total order $\leq$ on $J.$
For each finite subset $F \subset J,$
set $K_F = \prod_{i \in F} K (\QQ, n_i).$
We define a map $\af_F \colon K_F \to B_{\QQ}$
as follows.
Write $F = \{ i_1, \ldots, i_r \}$
with $i_1 < \cdots < i_r.$
Then, using $\cdot$ for the multiplication in $E_{\QQ},$ set
\[
\af_F (x_1, \ldots, x_r)
 = p_{\QQ} \left( \left( (\gm_{i_1} (x_1) \cdot \gm_{i_2} (x_2))
       \cdot \gamma_{i_3} (x_3) \right ) \cdot
        \cdots \cdot \gamma_{i_r} (x_r) \right).
\]
The induced homomorphism on homotopy groups restricts to the
inclusion $V_i \to \pi_* (B_{\QQ})$ on each summand of
$\pi_* (K_F) \cong \prod_{i \in F} V_i.$

If $F_1 \subset F_2 \subset J$ are finite sets,
there is an obvious map $k_{F_1, F_2} \colon K_{F_1} \to K_{F_2},$
obtained by viewing $K_{F_1}$ as the set of elements $x \in K_{F_2}$
such that $x_i$ is the basepoint of $K (\QQ, n_i)$
for $i \in F_2 \setminus F_1.$
Since the basepoints are the identities and the identities are
strict,
one checks that $\af_{F_2} \circ k_{F_1, F_2} = \af_{F_1}$
for $F_1 \subset F_2 \subset J.$
Accordingly, we obtain a map
$\alpha \colon \Dirlim_{F \subset J} K_F \to B_{\QQ},$
which moreover
induces an isomorphism on homotopy groups and hence is a homotopy
equivalence.
The H-space structures on the Eilenberg-Mac Lane spaces
give obvious H-space structures on the products $K_F,$
and the maps $k_{F_1, F_2}$ respect these structures
and the associated homotopies.
So $\Dirlim_{F \subset J} K_F$ is an H-space,
and thus $B_{\QQ}$ is also an H-space.
\end{proof}

\begin{Thm}
Let $A$ be commutative unital Banach algebra.
Let $\widetilde{\Lc_n (A)}\idc$ denote the universal cover of
$\Lc_n (A)\idc.$
If $\MM (A)$ is a finite complex then $\Lc_n (A)\idc$ is a rational
H-space.
In general, $\widetilde{\Lc_n (A)}\idc$ is a rational H-space.
\end{Thm}

\begin{proof}
When $\MM (A)$ is a finite complex the result is a consequence
of Theorem~\ref{HMRlocal}, Example~\ref{odd-sphere} and
Theorem~\ref{null-maps}, as argued before the statement
of Lemma~\ref{ratsur}.
For the general case, observe that,
by Theorem~\ref{T:les1} and Corollary~\ref{GL-injection},
the map $\widetilde{\GL_n (A)}\idc \to \widetilde{\Lc_n (A)}\idc$
induced on universal covers
by the natural
map $\GL_n (A)\idc \to \Lc_n (A)\idc$
is a surjection on rational homotopy groups.
As $\widetilde{\GL_n (A)}\idc$ is a topological group,
the result follows from Lemma~\ref{ratsur}.
\end{proof}

\end{document}